\documentclass[final]{siamltex}

\usepackage{graphicx}
\usepackage{amsmath, amssymb}
\usepackage[vlined,linesnumbered,algo2e]{algorithm2e}
\usepackage[usenames,dvipsnames,table]{xcolor}
\usepackage{booktabs}

\usepackage{tikz}
\usetikzlibrary{positioning,patterns,calc,arrows}

\usepackage{pgfplots,pgfplotstable}

\newcommand{\exptonumtickformat}{
  \pgfmathfloatparsenumber{\tick}
  \pgfmathfloatexp{\pgfmathresult}
  \pgfmathprintnumber{\pgfmathresult}
}

\usepackage{showkeys}              
\usepackage{lineno}                
\AtBeginDocument{%
\let\oldequation\equation%
\let\endoldequation\endequation%
\renewenvironment{equation}%
{\linenomath\oldequation}{\endoldequation\endlinenomath}%
\expandafter\let\expandafter\oldequationstar\csname equation*\endcsname%
\expandafter\let\expandafter\endoldequationstar\csname endequation*\endcsname%
\renewenvironment{equation*}%
{\linenomath\oldequationstar}{\endoldequationstar\endlinenomath}%
\let\oldalign\align%
\let\endoldalign\endalign%
\renewenvironment{align}%
{\linenomath\oldalign}{\endoldalign\endlinenomath}%
\expandafter\let\expandafter\oldalignstar\csname align*\endcsname%
\expandafter\let\expandafter\endoldalignstar\csname endalign*\endcsname
\renewenvironment{align*}%
{\linenomath\oldalignstar}{\endoldalignstar\endlinenomath}%
\let\oldflalign\flalign%
\let\endoldflalign\endflalign%
\renewenvironment{flalign}%
{\linenomath\oldflalign}{\endoldflalign\endlinenomath}%
\expandafter\let\expandafter\oldflalignstar\csname flalign*\endcsname%
\expandafter\let\expandafter\endoldflalignstar\csname endflalign*\endcsname%
\renewenvironment{flalign*}%
{\linenomath\oldflalignstar}{\endoldflalignstar\endlinenomath}%
\let\oldalignat\alignat%
\let\endoldalignat\endalignat%
\renewenvironment{alignat}%
{\linenomath\oldalignat}{\endoldalignat\endlinenomath}%
\expandafter\let\expandafter\oldalignatstar\csname alignat*\endcsname%
\expandafter\let\expandafter\endoldalignatstar\csname endalignat*\endcsname%
\renewenvironment{alignat*}%
{\linenomath\oldalignatstar}{\endoldalignatstar\endlinenomath}%
\let\oldgather\gather%
\let\endoldgather\endgather%
\renewenvironment{gather}%
{\linenomath\oldgather}{\endoldgather\endlinenomath}%
\expandafter\let\expandafter\oldgatherstar\csname gather*\endcsname%
\expandafter\let\expandafter\endoldgatherstar\csname endgather*\endcsname%
\renewenvironment{gather*}%
{\linenomath\oldgatherstar}{\endoldgatherstar\endlinenomath}%
\let\oldmultline\multline%
\let\endoldmultline\endmultline%
\renewenvironment{multline}%
{\linenomath\oldmultline}{\endoldmultline\endlinenomath}%
\expandafter\let\expandafter\oldmultlinestar\csname multline*\endcsname%
\expandafter\let\expandafter\endoldmultlinestar\csname endmultline*\endcsname%
\renewenvironment{multline*}%
{\linenomath\oldmultlinestar}{\endoldmultlinestar\endlinenomath}%
}                  
\usepackage[normalem]{ulem}	       

\usepackage[para,multiple]{footmisc}

\newcommand{\R}{\mathbb{R}}

\renewcommand{\O}{\mathcal{O}}
\newcommand{\abs}[1]{\left| {#1} \right|}

\renewcommand{\v}{\boldsymbol}
\newcommand{\vhat}[1]{\v{\hat{#1}}}
\renewcommand{\bar}{\overline}

\newcommand{\inner}[2]{#1 \cdot #2}
\renewcommand{\i}{\imath}

\renewcommand{\tilde}{\widetilde}
\newcommand{\conj}{\overline}
\renewcommand{\epsilon}{\varepsilon}
\newcommand{\eps}{\epsilon}
\renewcommand{\d}[1]{\,\mathrm{d} #1}

\usepackage{suffix}

\newcommand\secref[1]{Section~\ref{sec:#1}}
\WithSuffix\newcommand\secref*[1]{\ref{sec:#1}}

\newcommand\appref[1]{Appendix~\ref{app:#1}}
\WithSuffix\newcommand\appref*[1]{\ref{app:#1}}

\newcommand\tabref[1]{Table~\ref{tab:#1}}
\WithSuffix\newcommand\tabref*[1]{\ref{tab:#1}}

\newcommand\figref[1]{Figure~\ref{fig:#1}}
\WithSuffix\newcommand\figref*[1]{\ref{fig:#1}}

\newcommand\algref[1]{Algorithm~\ref{alg:#1}}
\WithSuffix\newcommand\algref*[1]{\ref{alg:#1}}

\newcommand\eqnref[1]{Equation~\eqref{eqn:#1}}
\WithSuffix\newcommand\eqnref*[1]{\eqref{eqn:#1}}

\newcommand{\F}{\mathcal{F}}   
\newcommand{\s}{\mathfrak{s}}  
\newcommand{\Ts}{T^{\s}}
\newcommand{\Tsl}{T^{\s,L}}
\newcommand{\Tsll}{T^{\s,LL}}
\newcommand{\FTs}{\tilde{T}^{\s}}
\newcommand{\eiksr}{e^{\i \kappa \inner{\vhat{s}}{\v{r}}}}   

\title{Fourier Based Fast Multipole Method for the Helmholtz Equation}

\author{Cris Cecka\thanks{Institute for Computational and Mathematical Engineering, Stanford University ({\tt }).} \and Eric Darve\thanks{Institute for Computational and Mathematical Engineering, Mechanical Engineering Department, Stanford University ({\tt darve@stanford.edu})}}
	
\begin{document}
\linenumbers

\maketitle

\begin{abstract}
    The fast multipole method (FMM) has had great success in reducing the computational complexity of solving the boundary integral form of the Helmholtz equation. We present a formulation of the Helmholtz FMM that uses Fourier basis functions rather than spherical harmonics. By modifying the transfer function in the precomputation stage of the FMM, time-critical stages of the algorithm are accelerated by causing the interpolation operators to become straightforward applications of fast Fourier transforms, retaining the diagonality of the transfer function, and providing a simplified error analysis. Using Fourier analysis, constructive algorithms are derived to a priori determine an integration quadrature for a given error tolerance. Sharp error bounds are derived and verified numerically. Various optimizations are considered to reduce the number of quadrature points and reduce the cost of computing the transfer function.
  \end{abstract}


\pagestyle{myheadings}
\thispagestyle{plain}
\markboth{C. CECKA AND E. DARVE}{FOURIER BASED FAST MULTIPOLE METHOD FOR THE HELMHOLTZ EQUATION}

\section{Introduction}

Since the development of the fast multipole method (FMM) for the wave equation in~\cite{Rokhlin_FMMScattering, Rokhlin_FMMPedestrian, Rokhlin_DiagonalTranslation, Rokhlin_FMM_EM, Rahola_DiagonalForms}, the FMM has proven to be a very effective tool for solving scalar acoustic and vector electromagnetic problems. In this paper, we consider the application of the FMM to the scalar Helmholtz equation, although our results can be immediately extended to the vector case as described in~\cite{Chew_Book,Darve_Numerical}. The application of the boundary element method to solve the integral form of the Helmholtz equation results in a dense linear system which can be solved by iterative methods such as GMRES or BCGSTAB. These methods require computing dense matrix-vector products which, using a direct implementation, are performed in $\O(N^2)$ floating-point operations. The FMM uses an approximation of the dense matrix to perform the product in $\O(N \log N)$ operations. This approximation is constructed from close-pair interactions and far-field approximations represented by spherical integrals that are accumulated and distributed through the domain via an octree.

There are a number of difficulties in implementing the FMM, each of which must be carefully considered and optimized to achieve the improved complexity. The most significant complication in the Helmholtz FMM is that the quadrature sampling rate must increase with the size of the box in the octree, requiring interpolation and anterpolation algorithms to transform the data between spherical quadratures of different levels of the tree. Local algorithms such as Lagrange interpolation and techniques which sparsify interpolant matrices are fast, but incur significant errors~\cite{Chew_ErrorDiagonal, Darve_Numerical}. Spherical harmonic transforms are global interpolation schemes and are exact but require fast versions for efficiency of the FMM. Many of these fast spherical transform algorithms are only approximate, complicated to implement or use, and not always stable~\cite{DriscollHealy_FST, HealyMoore_FST, Suda_FST}.

\begin{table}[t]
  \centering
  \begin{tabular}{cp{300pt}}
    \toprule
  Notation & Description\\\midrule
  $\kappa$ & wavenumber, $2\pi/\lambda$ with wavelength $\lambda$\\
  $a_l$ & box size at level $l$. Root: $l=0$. Highest active level: $l=2$.\\
  $\theta$ & polar angle \\
  $\phi$ & azimuthal angle \\
  $\conj{c}$ & complex conjugate of $c$\\
  $\v{x}$ & vector in $\R^3$, $\v{x} = \abs{\v{x}} \vhat{x}$\\
  $\inner{\v{x}}{\v{y}}$ & inner product, $\inner{\v{x}}{\v{y}} = \abs{\v{x}}\abs{\v{y}}\cos(\varphi_{\v{x},\v{y}})$\\
  $S^2$ & sphere, $\{ \vhat{s} \in \R^3 \, : \, \abs{\vhat{s}} = 1 \}$\\
  $j_n$ & spherical Bessel function of the first kind\\
  $y_n$ & spherical Bessel function of the second kind\\
  $h^{(1)}_n$ & spherical Hankel function of the first kind\\
  $P_n$ & Legendre polynomial\\
  $\F^x_n[f]$ & $n$th coefficient of $f$'s Fourier series in $x$, $f(x) = \sum_n \F^x_n[f] \, e^{\i n x}$\\
  \bottomrule
  \end{tabular}
\caption{Table of notations}
\label{tab:Notation}
\end{table}

In this paper, we use a multipole expansion which allows the use of 2D fast Fourier transforms (FFT) in the spherical coordinate system $(\phi,\theta)$. The main advantages are two fold: i) high performance FFT libraries are available on practically all computer platforms, resulting in accurate, robust, and fast interpolation algorithms; ii) the resulting error analysis is simplified and leads to sharp, a priori error bounds on the FMM. One of the difficulties in using FFTs is that we are forced to use a uniform distribution of points along $\phi$ and $\theta$ in the spherical quadrature. Naively, this leads to a much increased quadrature size for a given accuracy compared to the original spherical harmonics-based FMM. The reason is as follows. The multipole expansion in the high frequency regime is derived from:
\begin{equation*}
  \frac{e^{\i \kappa \abs{\v{r} + \v{r}_0}}}{\abs{\v{r} + \v{r}_0}} =
  \int_{\phi=0}^{2\pi} \int_{\theta=0}^{\pi} e^{\i\kappa\inner{\vhat{s}}{\v{r}}} \,
  T_{\ell,\v{r}_0}(\vhat{s}) \, \sin(\theta) \d{\theta}\d{\phi}
\end{equation*}
where $\vhat{s} = [\cos(\phi)\sin(\theta),\sin(\phi)\sin(\theta),\cos(\theta)]$ is the spherical unit vector. It is apparent that we are integrating along $\theta$ a function which has period $2\pi$. However the bounds of the integral are $0$ to $\pi$, over which interval the function has a discontinuity in its derivative. This results in a slow decay of the Fourier spectrum (essentially 1/freq$^2$) of the integrand. Consequently, a large number of quadrature points along $\theta$ are required.

We propose to use a variant of the scheme by J.\ Sarvas in~\cite{Sarvas_FFT} whereby the integration is extended from $0$ to $2\pi$ and the integrand modified:
\begin{equation}
  \frac{e^{\i \kappa \abs{\v{r} + \v{r}_0}}}{\abs{\v{r} + \v{r}_0}} =
  \frac{1}{2} \int_{\phi=0}^{2\pi} \int_{\theta=0}^{2\pi} e^{\i\kappa\inner{\vhat{s}}{\v{r}}} \,
  T_{\ell,\v{r}_0}(\vhat{s}) \, \abs{\sin(\theta)} \d{\theta} \d{\phi}
\label{eqn:modifiedIntegrand}
\end{equation}
We will describe in more details how an efficient scheme can be derived from this equation. The key property is that $e^{\i\kappa\inner{\vhat{s}}{\v{r}}}$ is approximately bandlimited in $\theta$ and therefore it is possible to remove the high frequency components of $T_{\ell,\v{r}_0}(\vhat{s}) |\sin(\theta)|$ without affecting the accuracy of the approximation. Using this smooth transfer function, which is now bandlimited in Fourier space, the number of quadrature points can be reduced dramatically. We show that the resulting number of quadrature points is reduced by about 40\% compared to the original spherical harmonics-based FMM. Consequently, we now have a scheme which requires few quadrature points and enables the use of efficient FFT routines. 

The approach in~\cite{Sarvas_FFT} is similar. However, rather than smoothing $T_{\ell,\v{r}_0}(\vhat{s}) \abs{\sin(\theta)}$ once during the precomputation phase as we detail in this paper, Sarvas instead incorporates the $\abs{\sin(\theta)}$ factor during the run-time phase of the FMM after the application of the transfer function. Although a detailed analysis is required for accurately assessing the relative cost of the two approaches, the technique in~\cite{Sarvas_FFT} requires approximately 1.5 times more sample points in the time-critical transfer pass of the algorithm, and requires an extra anterpolation step after applying $T$ with about twice more sample points than used for the method in this paper. We also note that the error analysis for the two methods is different, and is easier to carry out with the approach in this paper.

We derive a new a priori error analysis which incorporates both effects from truncation of the Gegenbauer series (a problem well analyzed~\cite{Chew_Book}) and the numerical quadrature. Our algorithm to predict the error is very sharp. The sharp bounds allow the method to choose a minimal number of quadrature points to guarantee a prescribed error. By comparison, the conventional approach leads to less accurate estimates resulting in either lower accuracy than requested or higher computational cost (over-estimation of the required approximation order). Although not considered in this paper, our error analysis approach can also be applied to the spherical harmonics-based FMM to yield similarly accurate error bounds. This has practical importance since it allows guaranteeing the error in the calculation while reducing the computational cost.

The novel contributions of this paper can be summarized as follows:
\begin{itemize}
\item Development of an efficient Helmholtz multi-level FMM which uses FFTs in the inter/anterpolation steps while retaining diagonal transfer and translation functions. The use of FFTs allows leveraging high performance FFT libraries available for most machines, sequential and parallel.
\item An error analysis that accounts for all error in the method and yields constructive algorithms to choose optimal method parameters.
\item Details of various optimizations to reduce the computational cost (e.g. use of symmetries in the precomputation of the transfer functions, use of symmetries on the unit sphere for the inter/anterpolation steps, optimization of the quadrature points near the poles of the unit sphere).
\item Pseudocodes are provided to clarify the method and help with an implementation by the reader.
\item Demonstration of the sharpness of the error bound and the asymptotic computational cost.
\end{itemize}

The paper is organized as follows. In \secref{HelmholtzFMM}, we introduce the critical parts of the classical Helmholtz FMM including the Gegenbauer series truncation (\secref*{GegenbauerTruncation}), the spherical quadrature (\secref*{SphereQuad}), and a short overview of interpolation/anterpolation strategies (\secref*{Interp}). 
\secref{FFMM} details the Fourier basis approach. The transfer function must be modified to lower the computational cost and obtain a competitive scheme, as detailed in \secref{MTF}. \secref{DirectQuadrature} analyzes the integration error to derive an algorithm which determines a quadrature with a prescribed error tolerance. The FFT based interpolation and anterpolation algorithms are described in \secref{FFTInterp} and numerical results are given in \secref{Results}. \tabref{Notation} lists the notations used in this paper.

\section{The Multilevel Fast Multipole Method}
\label{sec:HelmholtzFMM}

The FMM reduces the computational complexity of the matrix-vector multiplication
\begin{align}
\sigma_i = \sum_{j \neq i} \frac{e^{\i\kappa \abs{\v{x}_i -
      \v{x}_j}}}{\abs{\v{x}_i - \v{x}_j}} \psi_j
\label{eqn:HelmholtzMatrixVector}
\end{align}
for $i,j = 1, \ldots, N$ from $\O(N^2)$ to $\O(N\log N)$. This improvement is based on the Gegenbauer series
\begin{align}
\frac{e^{\i \kappa \abs{\v{r} + \v{r}_0}}}{\abs{\v{r} + \v{r}_0}} = \i \kappa \sum_{n=0}^\infty (-1)^n (2n+1) h_n^{(1)}(\kappa \abs{\v{r}_0})j_n(\kappa \abs{\v{r}})P_n(\inner{\vhat{r}}{\vhat{r}_0})
\label{eqn:Gegenbauer}
\end{align}
The series converges absolutely and uniformly for $\abs{\v{r}_0} \geq \frac{2}{\sqrt{3}} \abs{\v{r}}$ and has been studied extensively in~\cite{Collino_GS,Darve_ErrorAsymptotic}.

Truncating the Gegenbauer series at $\ell$
and using an integral over the unit sphere, $S^2$:
\begin{align*}
\frac{e^{\i \kappa \abs{\v{r} + \v{r}_0}}}{\abs{\v{r} + \v{r}_0}} = \int_{S^2}
e^{\i\kappa\inner{\vhat{s}}{\v{r}}} \, T_{\ell,\v{r}_0}(\vhat{s}) \d{S(\vhat{s})} + \epsilon_G
\end{align*}
where $\epsilon_G$ is the Gegenbauer series truncation error and the transfer function, $T_{\ell,\v{r}_0}(\vhat{s})$, is defined as
\begin{align}
T_{\ell,\v{r}_0}(\vhat{s}) = \frac{\i \kappa}{4\pi} \sum_{n=0}^\ell \i^n(2n+1)h_n^{(1)}(\kappa\abs{\v{r}_0})P_n(\inner{\vhat{s}}{\vhat{r}_0}).
\label{eqn:TransferFn}
\end{align}

The reduced computational complexity of the FMM is achieved by constructing a tree of nodes, typically an octree, over the domain of the source and field points. We recall the main steps of the FMM to set some notations. Let $M^l_{\alpha}(\vhat{s})$ be the outgoing field for $B^l_\alpha$, the box $\alpha$ of the tree in level $l \in [0,L]$ with center $\v{c}^l_\alpha$.\\
{\bf Initialization:} The method is initialized by computing the outgoing plane-wave expansions for each cluster contained in a leaf of the tree:
\begin{align*}
M^L_\alpha(\vhat{s}) = \sum_{i, \: \v{x}_i \in B^L_\alpha} \psi_i \, e^{\i \kappa \inner{\vhat{s}}{(\v{x}_i - \v{c}^L_\alpha)}}
\end{align*}
{\bf Upward Pass (M2M):} These outgoing expansions are then aggregated upward through the tree by accumulating the product of the child cluster expansions with the plane-wave translation function:
\begin{align*}
M^{l-1}_\alpha(\vhat{s}) = \sum_{\beta, \: B^{l}_\beta \subset B^{l-1}_\alpha} M^{l}_\beta(\vhat{s}) \, e^{\i \kappa \inner{\vhat{s}}{(\v{c}^{l}_\beta - \v{c}^{l-1}_\alpha)}} && l = L,\,L-1,\,\ldots,\,3
\end{align*}
{\bf Transfer Pass (M2L):} Incoming expansions, $I^l_\alpha(\vhat{s})$ of box $B^l_\alpha$, are computed from the outgoing by multiplication with the transfer function:
\begin{align*}
I^l_\alpha(\vhat{s}) = \sum_{\beta \in \mathcal{I}(B^l_\alpha)} M^l_\beta(\vhat{s}) \, T_{\ell,\v{c}^l_\beta - \v{c}^l_\alpha}(\vhat{s}) && l = L,\,L-1,\,\ldots,\,2
\end{align*}
where $\mathcal{I}(B^{l}_\alpha)$ is the interaction list of box $B^l_\alpha$, defined as all boxes of level $l$ which are not neighbors of $B^l_\alpha$, but whose parent is a neighbor of the parent of $B^l_\alpha$.\\
{\bf Downward Pass (L2L):} The incoming expansions are then disaggregated downward through the tree to compute the local field $L^l_\alpha(\vhat{s})$ of box $B^l_\alpha$:
\begin{align*}
L^{l+1}_\alpha(\vhat{s}) = L^{l}_\beta(\vhat{s}) \, e^{\i \kappa \inner{\vhat{s}}{(\v{c}^{l}_\beta - \v{c}^{l+1}_\alpha)}} + I^{l+1}_\alpha(\vhat{s}) && l = 2,\,3,\,\ldots,\,L-1
\end{align*}
where $B^{l+1}_\alpha \subset B^{l}_\beta$.\\
{\bf Field Computation:} At the finest level, the integration over the sphere is finally performed and added to the near-field contribution to determine the field value at the $N$ field points:
\begin{align}
\sigma_i = \int_{S^2} L^L_\alpha(\vhat{s})\, e^{\i \kappa \inner{\vhat{s}}{(\v{c}^L_\alpha - \v{x}_i)}} \d{S(\vhat{s})} \ \ + \sum_{\substack{j,\: j \neq i,\\ \v{x}_j \in \mathcal{N}(B^L_\alpha)}} \frac{e^{\i\kappa \abs{\v{x}_i - \v{x}_j}}}{\abs{\v{x}_i - \v{x}_j}} \psi_j
\label{eqn:FMM_FinalStep}
\end{align}
where $\v{x}_i \in B^L_\alpha$ and $\mathcal{N}(B^L_\alpha)$ is the neighbor list of $B^L_\alpha$, defined as $B^L_\alpha$ and all neighbor boxes of $B^L_\alpha$.

\subsection{Truncation Parameter in the FMM}
\label{sec:GegenbauerTruncation}

The truncation parameter $\ell$ must be chosen so that the Gegenbauer series \eqnref*{Gegenbauer} is converged to a desired accuracy. However, for $n > x$, $j_n(x)$ decreases super-exponentially while $h^{(1)}_n(x)$ diverges. The divergence of the Hankel function causes the transfer function to oscillate wildly and become numerically unstable. Even though the expansion converges, roundoff errors will adversely affect the accuracy if $\ell$ is too large. Thus, while one must choose $\ell > \kappa \abs{\v{r}}$ so that sufficient convergence of the Gegenbauer series is achieved, it must also be small enough to avoid the divergence of the transfer function. The selection of the truncation parameter $\ell$ has been studied extensively and a number of procedures for selecting it have been proposed~\cite{Rokhlin_FMMPedestrian,Darve_ErrorAsymptotic}.


The excess bandwidth formula (EBF) is derived from the convergence of the plane-wave spectrum as presented in~\cite{Chew_Book}. The EBF chooses $\ell$ as
\begin{align}
\ell \approx \kappa\abs{\v{r}} + C(\kappa\abs{\v{r}})^{1/3}
\label{eqn:EBF}
\end{align}
An empirically determined common choice is $C = 1.8 (d_0)^{2/3}$, where $d_0$ is the desired number of digits of accuracy. The EBF is one of the most popular choices to select the truncation parameter~\cite{Chew_ErrorControl}.


The actual Gegenbauer truncation error for a given $\ell$ can also be approximated.
As Carayol and Collino showed in~\cite{Collino_GS}, an upper bound of this error for large values of $\abs{\v{r}}$ is obtained when 
$P_n(\inner{\vhat{r}}{\vhat{r}_0}) = (\pm 1)^n$ 
so that
\begin{align}
\abs{\epsilon_G} &\lesssim \kappa \abs{\sum_{n=\ell+1}^\infty (\mp 1)^n (2n+1) h^{(1)}_n(\kappa\abs{\v{r}_0}) j_n(\kappa \abs{\v{r}})} \notag
\intertext{which they showed can be computed in closed form}
&= \kappa^2 \, \frac{\abs{\v{r}}\abs{\v{r}_0}}{\abs{\v{r}_0} \pm
  \abs{\v{r}}} \, \abs{h^{(1)}_{\ell+1}(\kappa\abs{\v{r}_0})j_{\ell}(\kappa\abs{\v{r}}) \pm h^{(1)}_{\ell}(\kappa\abs{\v{r}_0})j_{\ell+1}(\kappa\abs{\v{r}})}
\label{eqn:CollinoError}
\end{align}
This fails for small $\abs{\v{r}}$ when the upper bound is obtained by choosing $\inner{\vhat{r}}{\vhat{r}_0}$ such that the oscillation of $P_n(\inner{\vhat{r}}{\vhat{r}_0})$ compensates for the oscillation of $(-1)^n h^{(1)}_n(\kappa\abs{\v{r}_0}) j_n(\kappa \abs{\v{r}})$. Using the EBF as an initial guess for $\ell$ and refining the choice using the above closed form when $\abs{\v{r}}$ is sufficiently large is a simple algorithm which yields a nearly optimal value for $\ell$ (that is the smallest value consistent with the target error). This is the scheme we selected for this paper.

Carayol and Collino in~\cite{Collino_JA} and~\cite{Collino_GS} present an in-depth analysis of the Jacobi-Anger series and the Gegenbauer series. They find the asymptotic formula
\begin{align*}
  \ell \approx \kappa\abs{\v{r}} - \frac{1}{2} +
  \left(\frac{1}{2}\right)^{5/3}W^{2/3}\left( \frac{\kappa\abs{\v{r}}}{4\epsilon^6}
    \left(\frac{1+\abs{\v{r}_0}/\abs{\v{r}}}{1-\abs{\v{r}_0}/\abs{\v{r}}}\right)^{3/2}\right)
\end{align*}
where $W(x)$ is the Lambert function defined as the solution to
\begin{align*}
  W(x)e^{W(x)} = x \quad\quad x > 0
\end{align*}
This appears to be near optimal for large box sizes.

The errors introduced by this truncation have been investigated in other papers including~\cite{Chew_ErrorDiagonal,Collino_GS,Darve_ErrorAsymptotic}.

\subsection{Spherical Quadrature in the FMM}
\label{sec:SphereQuad}

The error analysis is simplified if a scheme is used which exactly integrates spherical harmonics, $Y^m_n$, up to some degree. 
The most common choice of quadrature uses uniform sample points in $\phi$ and Gauss-Legendre sample points in $z(\theta)$. With $N+1$ uniform points in the $\phi$ direction and $\frac{N+1}{2}$ Gauss-Legendre points in the $\theta$ direction, all $Y^m_n$, $-n \leq m \leq n$, $0 \leq n \leq N$ are integrated exactly~\cite{Darve_Numerical,Chew_ErrorDiagonal}.



\subsection{Interpolation and Anterpolation in the FMM}
\label{sec:Interp}

The quadrature sampling rate depends on the spectral content of the translation operator, $e^{\i \kappa \inner{\vhat{s}}{\v{r}}}$. Its coefficient in the spherical harmonic expansion decreases super-exponentially roughly for $n \gtrsim \kappa \abs{\v{r}}$. Therefore, as fields are aggregated in the upward pass and $\abs{\v{r}}$ becomes larger, a larger quadrature is required to resolve higher modes. These modes must be resolved since they interact with the modes in the transfer function, which do not significantly decay as $\ell$ increases.

Similarly, as fields are disaggregated in the downward pass, $\abs{\v{r}}$ becomes smaller and the higher modes of the incoming field make vanishingly small contributions to the integral as a consequence of Parseval's theorem. Thus, as the incoming field is disaggregated down the tree, a smaller quadrature can be used to resolve it. This makes the integration faster and is actually required to achieve an optimal asymptotic running time. See \tabref{complex}.

There have been several approaches to performing the interpolation and anterpolation between levels in the FMM. Below, we enumerate a number of options that have previously been studied. 

General, local interpolation methods like Lagrange interpolation, Gaussian interpolation, and B-splines are fast and provide for simple error analysis~\cite{Chew_ErrorDiagonal}.

A spherical harmonic transform maps function values $f_k$, sampled at $(\phi_k,\theta_k)$, to a new quadrature $(\phi'_{k'},\theta'_{k'})$, via the linear transformation
\begin{align}
f_{k'} = \sum_{m,l \leq K} Y^m_l(\phi'_{k'},\theta'_{k'}) \sum_k \omega_k \conj{Y^m_l(\phi_k,\theta_k)} f_k = \sum_k A_{k'k} f_k
\label{eqn:Interp}
\end{align}                    
This results in an $\O(N \log N)$ or $\O(N^{3/2})$ FMM (see \tabref{complex}). Fast spherical transforms (FST) have been developed in~\cite{DriscollHealy_FST,HealyMoore_FST,Suda_FST,RokhlinTygert_FST} and applied to the FMM in~\cite{Chowdhury_FST}. Using the FST reduces the interpolation and anterpolation procedures to $\O(K\log K)$, which results in an $\O(N \log^2 N)$ FMM. However, the accuracy and stability of these algorithms remain in question.

Approximations of the spherical transform have also been investigated in~\cite{AlpertJokobChien_SphereFilter,Darve_Numerical}. The interpolation matrix $A_{k'k}$ in \eqnref*{Interp} can be sparsified in a number of ways to provide an interpolation/anterpolation method that scales as $\O(K)$ with controllable relative error. Many other interpolation schemes exist with varying running times and errors. Rokhlin presents a fast polynomial interpolator based on the fast multipole method in~\cite{Rokhlin_FastPolynomialInterp}. See also~\cite{Knab_InterpProlateSpheroid}.

The asymptotic computational complexity in the big-O sense is summarized in \tabref{complex}.
            
\begin{table}[htbp]
	\centering
	\begin{tabular}{lll}
	Interpolation & Volume & Surface \\	\midrule
	Direct & $N \log N$ & $N^{3/2}$ \\
	Fast Global & $N$ & $N \log^2 N$ \\
	Local & $N$ & $N \log N$
	\end{tabular}
	\caption{Computational complexity in the big-O sense. Column 2 and 3 refer to the distribution of particles. We assume that $N = \O((\kappa a_0)^3)$ for a volume of scatterers and $N = \O((\kappa a_0)^2)$ for a surface of scatterers. Direct refers to a computation of the operator with no acceleration (basically several matrix vector products). Fast Global refers to fast spherical harmonic transforms or fast Fourier transforms (this work). Local are methods based on local interpolation such as using Lagrange interpolation. \label{tab:complex}}  
\end{table}

\section{Fourier Based Multilevel Fast Multipole Method}
\label{sec:FFMM}

The Fourier based fast multipole method is based on the identity
\begin{align}
\int_{S^2} \eiksr \, T_{\ell,\v{r}_0}(\vhat{s}) \d{S(\vhat{s})} = \int_0^{2\pi} \int_0^{2\pi} E_{\v{r}}(\theta,\phi) \, \Ts_{\ell,\v{r}_0}(\theta,\phi) \d{\phi} \d{\theta}
\label{eqn:FBFMM}
\end{align}
where the translation function $E_{\v{r}}(\theta,\phi)$ and the modified transfer function $\Ts_{\ell,\v{r}_0}(\theta,\phi)$ are defined as
\begin{align}
E_{\v{r}}(\theta,\phi) = \eiksr &&
\Ts_{\ell,\v{r}_0}(\theta,\phi) = \frac{1}{2}
\, T_{\ell,\v{r}_0}(\vhat{s})\abs{\sin(\theta)}
\end{align}
with $\vhat{s} = [\cos(\phi)\sin(\theta),\, \sin(\phi)\sin(\theta),\,\cos(\theta)]$ and $T_{\ell,\v{r}_0}(\vhat{s})$ is the transfer function defined in \eqnref{TransferFn}. The $\int_{S^2}$ representation is common and implied throughout the discussion in \secref{HelmholtzFMM}. The natural basis for integration on the sphere is the spherical harmonics $Y^m_n$ which form an orthonormal basis of $L^2(S^2)$ and the FMMs of \secref{HelmholtzFMM} attempt to preserve this basis expansion in the upward and downward pass. The $\int_0^{2\pi} \int_0^{2\pi}$ representation suggests the use of the Fourier functions $\{e^{\i n \phi} e^{\i m \theta}\}$ which form an orthonormal basis of $L^2([0,2\pi]\times[0,2\pi])$. Doubling the sphere to use Fourier methods is presented in a more general manner in~\cite{Sarvas_FFT,Sneeuw_2DFourierSphere}.

Using a Fourier basis rather than a spherical harmonics basis allows i) using two dimensional uniform quadratures; ii) fast Fourier transforms in the interpolation and anterpolation steps; and iii) spectral analysis in the error estimates. Of these advantages, the most important is that the FFT interpolations and anterpolations are fast and exact. Since there is no interpolation error, the only significant contributions to the final error are the truncation of the Gegenbauer series and the integration error due to the finite quadrature. Thus, the error analysis is simplified and we will determine in this paper precise bounds on the final error. In fact, our error analysis is fairly general and can be extended to the classical FMM with schemes that exactly integrate spherical harmonics (see direct and fast global methods in \secref{Interp}). The result is a fast, easy to implement, and controllable version of the FMM, which we detail in the following sections.

\subsection{Fourier Quadratures}
\label{sec:FourierBackground}

\subsubsection{Spherical Fourier}

A periodic complex function defined over $[0,2\pi] \times [0,2\pi]$ with spherical symmetry,
\begin{align}
f(\theta,\phi) = f(2\pi-\theta,\pi+\phi)
\label{eqn:RealSphereSym}
\end{align}
has the trigonometric polynomial representation
\begin{align}
f(\theta,\phi) = \sum_{n=-\infty}^{\infty} \sum_{m=-\infty}^{\infty} \tilde{f}_{n,m} e^{\i (n\theta+m\phi)}
\label{eqn:FourierPoly}
\end{align}
where the symmetry condition \eqnref*{RealSphereSym} is equivalent to
\begin{align}
\tilde{f}_{n,m} = (-1)^m \tilde{f}_{-n,m} \qquad \forall n,m
\label{eqn:FourierSphereSym}
\end{align}

If functions $f$ and $g$ can be represented exactly by trigonometric polynomials of degree $(\tilde{N}_\theta,\tilde{N}_\phi)$ then their $L^2$ inner-product can be computed exactly as
\begin{align*}
\langle f,g \rangle = \int_0^{2\pi} \int_0^{2\pi} f(\theta,\phi) \bar{g(\theta,\phi)} \d{\theta}\d{\phi} 
&= 4\pi^2 \sum_{n=-\tilde{N}_\theta}^{\tilde{N}_\theta} \sum_{m=-\tilde{N}_\phi}^{\tilde{N}_\phi} \tilde{f}_{n,m} \bar{\tilde{g}_{n,m}}\\
&= \frac{4\pi^2}{N_\theta N_\phi} \sum_{n=1}^{N_\theta} \sum_{m=1}^{N_\phi} f(\theta_n,\phi_m) \bar{g(\theta_n,\phi_m)}
\end{align*}
where 
\begin{align*}
N_\theta &= 2\tilde{N}_\theta + 1 & \theta_n &= \frac{2\pi n}{N_\theta}\\
N_\phi &= 2\tilde{N}_\phi + 1 & \phi_m &= \frac{2\pi m}{N_\phi}
\end{align*}
If $N_\phi$ is made even (by padding the Fourier coefficients with an extra zero), then
\begin{align}
f(\theta_n,\phi_m) = f(\theta_{N_\theta-n},\phi_{N_\phi/2 + m})
\label{eqn:NumSphereSym}
\end{align}
so that only half of the sampled values must be stored.

\subsubsection{Integration with Fourier Filtering}
\label{sec:IntFourierFilter}

A key concern in using the 2D Fourier basis rather than the spherical harmonics in \eqnref{FBFMM} is the integration weight $\frac{1}{2}\abs{\sin(\theta)}$. Although the transfer function, $T_{\ell,\v{r}_0}$ is bandlimited in $\theta$ and $\phi$ the modified transfer function, $\Ts_{\ell,\v{r}_0}$ is not bandlimited in $\theta$ due to the integration weight. In this section, we review a common strategy for integrating functions that are not bandlimited against functions that are bandlimited or nearly bandlimited.

Consider the periodic functions 
\begin{align*}
f(\theta) = \sum_{n = -\infty}^\infty \tilde{f}_n e^{\i n \theta} && 
g(\theta) = \sum_{m = -\infty}^\infty \tilde{g}_m e^{\i m \theta}
\end{align*}
Note that their exact integral on $[0,2\pi]$ is
\begin{align}
\langle f,g \rangle &= \int_0^{2\pi} f(\theta) \bar{g(\theta)} \d{\theta}
= 2\pi \sum_{n = -\infty}^\infty \tilde{f}_n \bar{\tilde{g}_n}
\label{eqn:Parseval}
\end{align}
Computing this integral numerically with a uniform quadrature of size $K$ yields
\begin{align*}
\langle f,g \rangle_K = \frac{2\pi}{K} \sum_{k = 1}^{K} f\left(\theta_k\right) \bar{g\left(\theta_k\right)}
&= \frac{2\pi}{K} \sum_{n = -\infty}^\infty \sum_{m = -\infty}^\infty \tilde{f}_n \bar{\tilde{g}_m} \sum_{k = 1}^{K} e^{2\pi\i(n-m)k/K}\\
&= 2\pi \sum_{n = -\infty}^\infty \sum_{m = -\infty}^\infty \tilde{f}_n \bar{\tilde{g}_{n + mK}}
\end{align*}
So the error is exactly given by
\begin{align*}
\abs{\langle f,g \rangle - \langle f,g \rangle_K} = 2\pi \abs{\sum_{n = -\infty}^\infty \sum_{m \neq 0} \tilde{f}_n \bar{\tilde{g}_{n + mK}}}
= 2\pi \abs{\sum_{n = -\infty}^\infty \sum_{m \neq 0} \tilde{f}_{n + mK} \bar{\tilde{g}_n}}
\end{align*}
The error is therefore determined by the asymptotic decay of both $\tilde{f}_n$ and $\tilde{g}_n$. For example, if $|\tilde{f}_n|$ decays as $1/n^\alpha$, $\alpha > 1$, then $\sum_n |\tilde{f}_{nK} \tilde{g}_0| \sim 1/K^\alpha$, which is the typically expected decay of the error on a uniform quadrature as a function of $K$. The equivalent result holds for $\tilde{g}_n$. Thus, if the spectrum of either $f$ or $g$ decays slowly, a very large quadrature $K$ is needed.

However, if the spectrum of $f$ decays slowly and the spectrum of $g$ decays quickly a much smaller quadrature can be used by truncating the Fourier coefficients of $f$. Let us define a bandlimited version of $f(\theta)$
\begin{equation}
f_N(\theta) = \sum_{n = -N}^N  \tilde{f}_n e^{\i n \theta}
\label{eqn:fN}
\end{equation}
and use a uniform quadrature of $K$ points to compute the integral. Then
\begin{equation*}
\langle f_N,g \rangle_K = 2\pi \sum_{n = -N}^N \sum_{m = -\infty}^\infty \tilde{f}_n \bar{\tilde{g}_{n + mK}}
\end{equation*}
which yields the error
\begin{align}
\abs{\langle f,g \rangle - \langle f_N,g \rangle_K} = 2\pi \abs{\sum_{\abs{n} > N} \tilde{f}_n \bar{\tilde{g}_n} - \sum_{n = -N}^N \sum_{m \neq 0} \tilde{f}_n \bar{\tilde{g}_{n + mK}}}
\label{eqn:FourierError}
\end{align}
This error can be made small by requiring only that $\tilde{g}_n$ decays quickly enough. That is, $|\tilde{f}_n|$ can decay slowly provided $|\tilde{g}_n|$ decays quickly. The first term of \eqnref*{FourierError} represents the truncation error in using the bandlimited $f_N$ instead of $f$. The second term represents the aliasing error resulting from the finite sampling of the function $g$. As an example, if we assume that $\tilde{g}_n$ is negligible for $\abs{n} > N$ then it is sufficient to choose $K = 2N+1$ so that $\abs{n + mK} > N$ for all $-N \leq n \leq N$ and $m \neq 0$ so that $\tilde{g}_{n + mK}$ is always negligible.

\begin{figure}[htbp]
\centering
%
%
\begin{tikzpicture}

\definecolor{mycolor1}{rgb}{0,0.5,0}

\begin{loglogaxis}[%
scale only axis,
width=4in,
height=1.5in,
grid=both,
xticklabel={\exptonumtickformat},
xmin=1, xmax=4096,
ymin=2e-16, ymax=10000,
legend entries={Fourier Filter, Trapezoid, $90/K^2$},
legend pos = south west,
xlabel=$N$,
ylabel=Relative Error]

\addplot [
color=blue,
solid,
mark=+,
mark options={solid} 
]
coordinates{
 (1,47.4111)(2,53.6121)(2,53.6121)(2,53.6121)(3,40.4771)(3,40.4771)(4,26.5787)(5,52.1291)(6,12.2683)(7,30.6811)(8,17.4813)(10,33.5354)(11,20.2775)(13,1.29091)(16,5.1609)(19,11.01)(23,16.9949)(27,8.1178)(32,14.3957)(38,0.0339386)(45,0.00394055)(54,0.000112673)(64,2.59461e-05)(76,7.48466e-08)(91,3.03109e-11)(108,6.25086e-14)(128,1.65322e-14)(152,1.12226e-14)(181,3.25817e-15)(215,2.4738e-14)(256,2.89615e-15)(304,8.10923e-14)(362,1.04986e-14)(431,4.22356e-15)(512,1.68942e-15)(609,9.89519e-15)(724,6.39567e-15)(861,1.02572e-14)(1024,4.21149e-14)(1218,1.20673e-16)(1448,3.01683e-15)(1722,4.58557e-15)(2048,5.18894e-15)(2435,1.8101e-15)(2896,1.20673e-16)(3444,4.4649e-15)(4096,1.20673e-16) 
};

\addplot [
color=mycolor1,
only marks,
mark=square,
mark options={solid},
mark size=2.5
]
coordinates{
 (1,53.6348)(2,51.6286)(2,51.6286)(2,51.6286)(3,37.6699)(3,37.6699)(4,25.4531)(5,51.4105)(6,12.6303)(7,29.717)(8,16.9676)(10,33.7142)(11,19.9617)(13,1.31468)(16,5.17657)(19,10.9612)(23,16.9543)(27,8.16758)(32,14.4126)(38,0.0490339)(45,0.0102191)(54,0.00722831)(64,0.00521278)(76,0.00373389)(91,0.00262557)(108,0.00187522)(128,0.00134124)(152,0.000954699)(181,0.000675323)(215,0.000479765)(256,0.000339057)(304,0.000240814)(362,0.000170047)(431,0.000120084)(512,8.51666e-05)(609,6.02393e-05)(724,4.26471e-05)(861,3.01695e-05)(1024,2.13377e-05)(1218,1.50868e-05)(1448,1.06776e-05)(1722,7.5517e-06)(2048,5.33992e-06)(2435,3.77804e-06)(2896,2.67132e-06)(3444,1.88906e-06)(4096,1.33565e-06) 
};

\addplot [
color=black,
dashed
]
coordinates{
 (1,5.625)(2,2.5)(2,2.5)(2,2.5)(3,1.40625)(3,1.40625)(4,0.9)(5,0.625)(6,0.459184)(7,0.351562)(8,0.277778)(10,0.18595)(11,0.15625)(13,0.114796)(16,0.0778547)(19,0.05625)(23,0.0390625)(27,0.028699)(32,0.0206612)(38,0.0147929)(45,0.0106333)(54,0.00743802)(64,0.00532544)(76,0.00379491)(91,0.00265832)(108,0.00189378)(128,0.00135208)(152,0.000961169)(181,0.000679266)(215,0.000482253)(256,0.000340656)(304,0.00024187)(362,0.000170753)(431,0.000120563)(512,8.54964e-05)(609,6.04676e-05)(724,4.28062e-05)(861,3.02808e-05)(1024,2.14158e-05)(1218,1.51417e-05)(1448,1.07163e-05)(1722,7.579e-06)(2048,5.35918e-06)(2435,3.79165e-06)(2896,2.68093e-06)(3444,1.89585e-06)(4096,1.34045e-06) 
};

\end{loglogaxis}
\end{tikzpicture}	
\caption{Relative error in the integral as a function of $N$. The number of quadrature points is $K=2N+1$. The blue line integrates the function using a Fourier filter for $|\sin(\theta)|$ as shown in \eqnref{FourierError}. The rapid convergence can be seen after $N \gtrsim 66$. The trapezoid method has a slower convergence in $1/K^2$ because $|\sin \theta|$ is only $C^0$.} 
\label{fig:SINconv}
\end{figure}
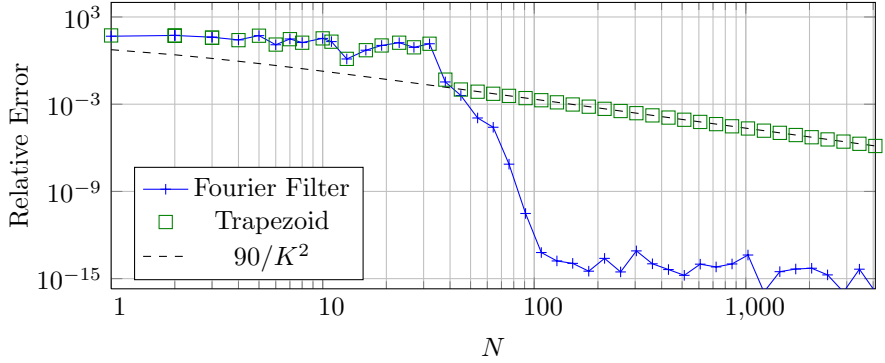

Detailed numerical results will be presented in Section~\ref{sec:Results}. However, to demonstrate this idea on a simple example, we build a quadrature to calculate:
\begin{equation*}
	\int_0^{2\pi} \abs{\sin \theta} \, \cos(64 \cos(\theta)) \, \d\theta
\end{equation*}
which is a simple model problem for [see~\eqnref{FBFMM}]
\begin{equation*}
	\int_0^{2\pi} T_{\ell,\v{r}_0}(\vhat{s})\abs{\sin \theta} \, \eiksr \, \d{\theta}
\end{equation*}  
when $\v{r} = 64/\kappa \; \vhat{z}$ ($z$-axis). Using a uniform quadrature to compute this integral yields slow convergence in $1/N^2$ because the function $f(\theta)=\abs{\sin \theta}$ is $C^0$ (see \figref{SINconv}). Instead, if we smooth $\abs{\sin \theta}$ and remove the high-frequency components following \eqnref{fN}, we obtain a much faster convergence (see \figref{SINconv}). The number of quadrature points is $K=2N+1$. The Fourier spectrum of $g(\theta)=\cos(64 \cos(\theta))$ decays rapidly once $\abs{n} \gtrsim 64$ (see \figref{SINFspec}). Convergence should occur once $N \gtrsim 64$ and this is indeed what is observed in~\figref{SINconv}. The exact value of the integral $\sin(64)/16$ is used as the reference solution to calculate the error.


\begin{figure}[htbp]
\centering
%
%
\begin{tikzpicture}

\begin{semilogyaxis}[%
scale only axis,
width=4in,
height=1.25in,
grid=both,
xmin=-128, xmax=128,
ytick={1e0,1e-4,1e-8,1e-12,1e-16},
ymin=1e-17, ymax=1e0,
xlabel=Fourier frequency $n$,
ylabel=$\F^\theta_n(g(\theta))$]

\addplot [
color=blue,
solid
]
coordinates{
 (-128,1.17961e-16)(-127,4.22686e-16)(-126,7.06082e-16)(-125,5.57722e-16)(-124,5.05725e-16)(-123,2.27037e-16)(-122,2.91535e-16)(-121,2.18467e-16)(-120,7.66082e-17)(-119,6.24088e-17)(-118,2.89562e-16)(-117,2.06891e-16)(-116,2.01499e-16)(-115,4.96582e-16)(-114,2.48108e-16)(-113,5.09118e-16)(-112,8.78985e-17)(-111,3.17371e-16)(-110,5.85495e-16)(-109,9.57796e-16)(-108,6.46762e-16)(-107,1.08032e-15)(-106,1.54355e-15)(-105,1.12427e-15)(-104,1.19104e-14)(-103,6.79472e-16)(-102,9.70319e-14)(-101,3.08421e-16)(-100,7.72111e-13)(-99,6.37196e-16)(-98,5.83869e-12)(-97,2.34591e-16)(-96,4.18857e-11)(-95,5.32373e-16)(-94,2.84417e-10)(-93,2.5816e-16)(-92,1.82425e-09)(-91,2.73661e-16)(-90,1.10271e-08)(-89,4.17343e-16)(-88,6.2661e-08)(-87,2.58911e-16)(-86,3.33795e-07)(-85,8.04502e-16)(-84,1.66172e-06)(-83,8.96233e-16)(-82,7.70369e-06)(-81,7.55527e-16)(-80,3.31252e-05)(-79,6.31724e-16)(-78,0.000131499)(-77,1.21256e-16)(-76,0.000479316)(-75,5.20502e-16)(-74,0.0015938)(-73,4.20897e-17)(-72,0.00479629)(-71,4.75068e-16)(-70,0.0129327)(-69,2.78051e-16)(-68,0.0308385)(-67,2.86259e-16)(-66,0.0638664)(-65,1.43644e-16)(-64,0.111821)(-63,3.13855e-16)(-62,0.158192)(-61,7.01642e-16)(-60,0.164628)(-59,3.68277e-16)(-58,0.0922594)(-57,3.10397e-16)(-56,0.0425775)(-55,4.46787e-16)(-54,0.133426)(-53,3.06251e-16)(-52,0.0698861)(-51,5.05406e-16)(-50,0.0845323)(-49,7.39837e-16)(-48,0.103646)(-47,5.87549e-16)(-46,0.0557984)(-45,2.67696e-16)(-44,0.102809)(-43,6.73895e-16)(-42,0.0644115)(-41,6.57048e-16)(-40,0.0805168)(-39,2.6799e-16)(-38,0.0957131)(-37,3.76125e-16)(-36,0.0212885)(-35,1.90587e-16)(-34,0.105771)(-33,1.09656e-16)(-32,0.0695322)(-31,5.95673e-16)(-30,0.0318695)(-29,2.70271e-16)(-28,0.0996525)(-27,2.65585e-16)(-26,0.0891893)(-25,6.04045e-16)(-24,0.0228868)(-23,2.42283e-16)(-22,0.0504489)(-21,6.03455e-16)(-20,0.0946464)(-19,4.79545e-16)(-18,0.099512)(-17,8.25749e-16)(-16,0.0741284)(-15,7.45418e-16)(-14,0.0343573)(-13,6.42841e-16)(-12,0.0062175)(-11,4.55328e-16)(-10,0.0397485)(-9,3.16398e-16)(-8,0.0636895)(-7,4.85728e-16)(-6,0.0788272)(-5,3.75814e-16)(-4,0.0873305)(-3,8.06647e-16)(-2,0.091409)(-1,4.80193e-16)(0,0.09259)(1,4.80193e-16)(2,0.091409)(3,8.06647e-16)(4,0.0873305)(5,3.75814e-16)(6,0.0788272)(7,4.85728e-16)(8,0.0636895)(9,3.16398e-16)(10,0.0397485)(11,4.55328e-16)(12,0.0062175)(13,6.42841e-16)(14,0.0343573)(15,7.45418e-16)(16,0.0741284)(17,8.25749e-16)(18,0.099512)(19,4.79545e-16)(20,0.0946464)(21,6.03455e-16)(22,0.0504489)(23,2.42283e-16)(24,0.0228868)(25,6.04045e-16)(26,0.0891893)(27,2.65585e-16)(28,0.0996525)(29,2.70271e-16)(30,0.0318695)(31,5.95673e-16)(32,0.0695322)(33,1.09656e-16)(34,0.105771)(35,1.90587e-16)(36,0.0212885)(37,3.76125e-16)(38,0.0957131)(39,2.6799e-16)(40,0.0805168)(41,6.57048e-16)(42,0.0644115)(43,6.73895e-16)(44,0.102809)(45,2.67696e-16)(46,0.0557984)(47,5.87549e-16)(48,0.103646)(49,7.39837e-16)(50,0.0845323)(51,5.05406e-16)(52,0.0698861)(53,3.06251e-16)(54,0.133426)(55,4.46787e-16)(56,0.0425775)(57,3.10397e-16)(58,0.0922594)(59,3.68277e-16)(60,0.164628)(61,7.01642e-16)(62,0.158192)(63,3.13855e-16)(64,0.111821)(65,1.43644e-16)(66,0.0638664)(67,2.86259e-16)(68,0.0308385)(69,2.78051e-16)(70,0.0129327)(71,4.75068e-16)(72,0.00479629)(73,4.20897e-17)(74,0.0015938)(75,5.20502e-16)(76,0.000479316)(77,1.21256e-16)(78,0.000131499)(79,6.31724e-16)(80,3.31252e-05)(81,7.55527e-16)(82,7.70369e-06)(83,8.96233e-16)(84,1.66172e-06)(85,8.04502e-16)(86,3.33795e-07)(87,2.58911e-16)(88,6.2661e-08)(89,4.17343e-16)(90,1.10271e-08)(91,2.73661e-16)(92,1.82425e-09)(93,2.5816e-16)(94,2.84417e-10)(95,5.32373e-16)(96,4.18857e-11)(97,2.34591e-16)(98,5.83868e-12)(99,6.37196e-16)(100,7.72115e-13)(101,3.08421e-16)(102,9.70317e-14)(103,6.79472e-16)(104,1.19142e-14)(105,1.12427e-15)(106,1.54374e-15)(107,1.08032e-15)(108,6.41572e-16)(109,9.57796e-16)(110,5.96991e-16)(111,3.17371e-16)(112,8.78985e-17)(113,5.09118e-16)(114,2.66703e-16)(115,4.96582e-16)(116,1.91988e-16)(117,2.06891e-16)(118,2.94595e-16)(119,6.24088e-17)(120,7.60978e-17)(121,2.18467e-16)(122,2.88595e-16)(123,2.27037e-16)(124,5.00446e-16)(125,5.57722e-16)(126,6.9643e-16)(127,4.22686e-16) 
};

\end{semilogyaxis}
\end{tikzpicture}	
\caption{The spectrum of $g(\theta) = \cos(64 \cos(\theta))$ showing a rapid decay of the coefficients for $|n| > 64$. The Fourier spectrum was computed using 256 sample points.}
\label{fig:SINFspec}
\end{figure}
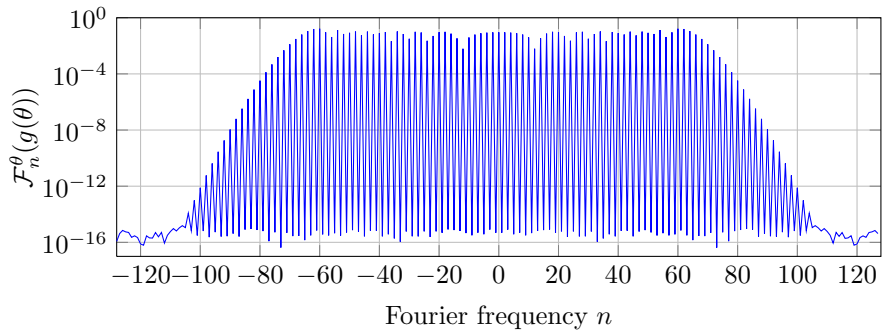

The idea of accurately calculating the integral of a product of two functions by analytically removing high frequencies in one of the two functions can be found in other papers dealing with the fast multipole method for the Helmholtz kernel in the high frequency regime, e.g.~\cite{Darve_Numerical,Darve_ErrorAsymptotic}. In the context of these papers, the smoothing operation (removal of high frequencies) is often termed anterpolation or subsampling. A similar idea is found in Sarvas et al.~\cite{Sarvas_FFT}. In McKay Hyde et al.~\cite{mckay2005fast} (Appendix A, p.~254--257), this idea is used in the more general context of calculating the integral of the product of a discontinuous function with a $C^1$ piecewise-smooth and periodic function.

\subsubsection{Fourier Interpolation and Anterpolation}

Fast Fourier interpolation and anterpolation methods are key to creating an efficient Fourier based FMM by supersampling trigonometric polynomials and truncating spectral content that does not significantly contribute to the final result. A Fourier interpolation pads the Fourier coefficients of a function with zeros and increases the sampling rate in real-space. A Fourier anterpolation removes high frequencies of a function and decreases the sampling rate in real-space.

To motivate the use of Fourier interpolations and anterpolations in the Helmholtz MLFMM, suppose that the spectrum of $f_N(\theta)$ is bounded, that is $|\tilde{f}_n| \leq F$, and that $K = 2N+1$, then \eqnref*{FourierError} simplifies to
\begin{align}
\abs{\langle f,g \rangle - \langle f_N,g \rangle_K} \leq 4\pi F \abs{\sum_{\abs{n} > N} \tilde{g}_n}.
\label{eqn:NBound}
\end{align}
This can be used to find an appropriate truncation parameter $N$ if $\tilde{g}_n$ is known or can be approximated.
The key step to constructing a fast algorithm is to note that the Fourier series of $E_{\v{r}}(\theta,\phi)$ in $\theta$ decays rapidly while the Fourier series of $\Ts_{\ell,\v{r}_0}(\theta,\phi)$ in $\theta$ decays slowly. This is due to the slow decay of the Fourier series of $\abs{\sin(\theta)}$:
\begin{align*}
\F^\theta_n[\abs{\sin(\theta)}] = \frac{(-1)^n + 1}{\pi(1 - n^2)} =
\begin{cases}
\frac{2}{\pi}\frac{1}{1-n^2} & \text{if $n$ even}\\
0 & \text{if $n$ odd}
\end{cases}
\end{align*}
The spectrum of the plane-wave is given by
\begin{align}
E_{\v{r}}(\theta,\phi) = e^{\i\kappa\abs{\v{r}}\cos(\varphi_{\vhat{s},\v{r}})} = \sum_{n=-\infty}^\infty \i^n J_n(\kappa\abs{\v{r}})e^{\i n \varphi_{\vhat{s},\v{r}}}
\label{eqn:PlaneWaveSpec}
\end{align}
where $\varphi_{\vhat{s},\v{r}}$ is the angle between $\vhat{s}$ and $\v{r}$. The functions $J_n$ decay rapidly once $n > \kappa\abs{\v{r}}$ resulting in an approximately band-limited function. In the following sections, we use the spectral decay of the translation function $E_{\v{r}}(\theta,\phi)$ (which plays the role of $g$ in the previous section) to determine an appropriate Fourier truncation of the modified transfer function, $\Ts_{\ell,\v{r_0}}(\theta,\phi)$ (which plays the role of $f$).

\eqnref{NBound} illustrates the need for interpolation and anterpolation as described in \secref{Interp} but in the context of the Fourier basis. With $\tilde{g}_n = \F_n^\theta[E_{\abs{\v{r}}\vhat{z}}] = \i^n J_n(\kappa\abs{\v{r}})$, the truncation $N$ is approximately $N \gtrsim \kappa\abs{\v{r}} \sim \ell$. During the upward pass, $\abs{\v{r}}$ increases and we require more modes in $f = \Ts_{\ell,\v{r}_0}$. During the downward pass, the incoming local field, $L^l$ in \secref{HelmholtzFMM}, is to be integrated against translation functions of increasingly smaller $\abs{\v{r}}$. See \eqnref{FMM_FinalStep} with $\v{r} = \v{c}^L_\alpha - \v{x}_i$. The high modes of the field do not significantly contribute to the exact integral \eqnref*{Parseval} so the field can be safely anterpolated at each downward step.

\subsection{Computing the Bandlimited Modified Transfer Function}
\label{sec:MTF}

Select a 2D uniform quadrature with points $(\theta_n,\phi_m)$ defined by
\begin{align*}
\theta_n = \frac{2\pi n}{N_\theta} \qquad \qquad \phi_m = \frac{2\pi m}{N_\phi}
\end{align*}

Noting that the plane wave $E_{\v{r}}(\theta,\phi)$ and modified transfer function $\Ts_{\ell,\v{r}_0}(\theta,\phi)$ both have spherical symmetry \eqnref*{RealSphereSym}, the computational and memory cost are reduced by requiring $N_\phi$ to be even so that only half of the quadrature points need to be computed and stored.

Additionally, in an FMM with a single octree, there are 316 distinct transfer vectors $\v{r}_0$ per level. By enforcing symmetries in the quadrature, the number of modified transfer functions that must be precomputed is reduced. Specifically, by requiring $N_\theta$ to be a multiple of 2 and $N_\phi$ to be a multiple of 4, we enforce reflection symmetries in the $z = 0$, $x = 0$, $y = 0$, $x = y$, and $x = -y$ planes. This reduces the number of modified transfer functions that need to be precomputed from 316 per level to 34 -- saving a factor of 9.3 in memory and costing a negligible permutation of the values of a computed modified transfer function. See \figref{SemiOctant}.

\setlength{\unitlength}{.4mm} 
\begin{figure}
\centering
	\begin{picture}(0, 120)(0,0) 
\color{black}
\thicklines
\put(-130,10){\vector(1,0){26}}
\put(-104,12){{\small x}}
\put(-130,10){\vector(-1,-1){13}}
\put(-135,-3){{\small y}}
\put(-130,10){\vector(0,1){26}}
\put(-128,36){{\small z}}
\put(-20, 40){\line(0,1){80}}
\put(-40, 40){\line(0,1){80}}
\put(-60, 40){\line(0,1){80}}
\put(-60, 40){\line(1,0){40}}
\put(-60, 60){\line(1,0){40}}
\put(-100,120){\line(-1,-1){10}}
\put(-100,100){\line(-1,-1){10}}
\put(-100,80){\line(-1,-1){10}}
\put(-80,80){\line(-1,-1){20}}
\put(-60,80){\line(-1,-1){30}}
\put(-60,60){\line(-1,-1){30}}
\put(-60,40){\line(-1,-1){30}}
\put(-20,40){\line(-1,-1){40}}
\put(-40,40){\line(-1,-1){40}}
\put(-100,60){\line(1,0){20}}
\put(-100,60){\line(0,1){10}}
\put(-90,50){\line(0,1){10}}
\put(-80,0){\line(0,1){10}}
\put(-90,50){\line(0,-1){40}}
\put(-110,110){\line(0,-1){40}}
\put(-80,20){\line(0,1){40}}
\put(-70,30){\line(0,1){40}}
\put(-70,30){\line(1,0){40}}
\put(-80,20){\line(1,0){40}}
\put(-90,10){\line(1,0){40}}
\put(-80,0){\line(1,0){20}}
\put(-100,80){\line(0,1){40}}
\put(-80,80){\line(0,1){40}}
\put(-100,120){\line(1,0){80}}
\put(-100,100){\line(1,0){80}}
\put(-100,80){\line(1,0){80}}
\put(-110,70){\line(1,0){40}}
%
\put(140,0){\vector(1,0){26}}
\put(165,2){{\small x}}
\put(140,0){\vector(1,1){13}}
\put(155,13){{\small y}}
\put(140,0){\vector(0,1){26}}
\put(142,25){{\small z}}
\put(60,0){\line(0,1){80}}
\put(80,0){\line(0,1){80}}
\put(100,0){\line(0,1){80}}
\put(110,10){\line(0,1){80}}
\put(120,20){\line(0,1){80}}
\put(130,30){\line(0,1){80}}
\put(140,40){\line(0,1){80}}
\put(100,0){\line(-1,0){40}}
\put(100,20){\line(-1,0){40}}
\put(100,40){\line(-1,0){80}}
\put(100,60){\line(-1,0){80}}
\put(100,80){\line(-1,0){80}}
\put(20,40){\line(0,1){40}}
\put(40,40){\line(0,1){40}}
\put(100,0){\line(1,1){40}}
\put(100,20){\line(1,1){40}}
\put(100,40){\line(1,1){40}}
\put(100,60){\line(1,1){40}}
\put(100,80){\line(1,1){40}}
\put(80,80){\line(1,1){40}}
\put(60,80){\line(1,1){30}}
\put(40,80){\line(1,1){20}}
\put(20,80){\line(1,1){10}}
\put(110,90){\line(-1,0){80}}
\put(120,100){\line(-1,0){60}}
\put(130,110){\line(-1,0){40}}
\put(140,120){\line(-1,0){20}}
\color{gray}
\put(20,0){\line(0,1){20}}
\put(40,0){\line(0,1){20}}
\put(20,0){\line(1,0){20}}
\put(20,20){\line(1,0){20}}
\put(20,20){\line(1,1){10}}
\put(40,20){\line(1,1){10}}
\put(40,0){\line(1,1){10}}
\put(50,10){\line(0,1){20}}
\put(30,30){\line(1,0){20}}
\end{picture}
\caption{The center of each box represents one transfer vector $\v{r}_0$ which must be computed. The pictures on the left and right panels represent the same set of boxes viewed under two different angles. Due to the symmetries of the quadrature, we need only compute transfer vectors with $x,y,z \geq 0$ and $x \geq y$. We therefore end up with essentially half of an octant. Specifically, 34 transfer vectors are required; they can be reflected into any of the 316 needed.}
\label{fig:SemiOctant}
\end{figure}
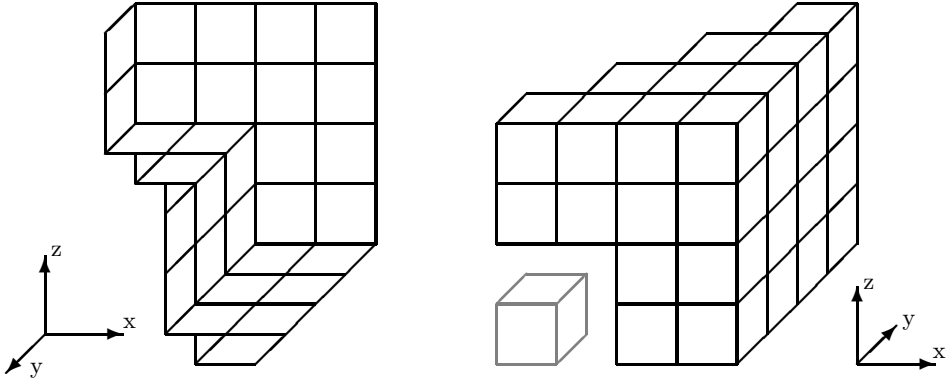

Suppose we have chosen a quadrature characterized by $(N_\theta,N_\phi)$. Following \secref{IntFourierFilter}, we need to exactly calculate a bandlimited version of $\Ts_{\ell,\v{r}_0}$, called $\Tsl_{\ell,\v{r}_0}$, such that:
\begin{equation*}
  \F^\theta_n[\Tsl_{\ell,\v{r}_0}(\theta,\phi))] =
  \begin{cases}
      \F^\theta_n[\Ts_{\ell,\v{r}_0}(\theta,\phi)], & \;\; \abs{n} \leq N_\theta/2-1 \\
      0, & \; \; \text{otherwise}
  \end{cases}
\end{equation*}
Since $T_{\ell,\v{r}_0}$ is bandlimited in $\theta$ with bandwidth $2\ell+1$, only the frequencies $\abs{m} \leq N_\theta/2+\ell-1$ of $\abs{\sin(\theta)}$ contribute to $\Tsl_{\ell,\v{r}_0}$. Therefore, the exact bandlimited modified transfer function, $\Tsl_{\ell,\v{r}_0}$, can be computed using the pseudocode in \algref{BMTF}.

\begin{algorithm2e}[H]
$\tilde{s}_n \leftarrow \F^\theta_n[\abs{\sin(\theta)}] = \frac{(-1)^n + 1}{\pi(1 - n^2)}$ for all $\abs{n} \leq N_\theta/2 + \ell - 1$\;
\For{$\phi_m$, $m = 0, \ldots, N_\phi/2-1$,} {
    $T(\theta_n,\phi_m) \leftarrow \frac{1}{2} T_{\ell,\v{r}_0}(\frac{2\pi n}{2\ell+1},\phi_m)$, $n = 0, \ldots, 2\ell$\;
    $\tilde{T}_n \leftarrow \F^\theta_n[T]$\;

    $\Tsl_n \leftarrow \tilde{s} \otimes \tilde{T}$ convolution of Fourier series\;
    $\Tsl_n \leftarrow$ truncate to frequencies $\abs{n} \leq N_\theta/2-1$\;
    $\Tsl(\theta_n,\phi_m) \leftarrow $ inverse transform of $\Tsl_n$\;
}
\caption{Pseudocode to compute the bandlimited modified transfer function, $\Tsl$, given $\ell$, $\v{r}_0$, $N_\theta$, and $N_\phi$.}
\label{alg:BMTF}
\end{algorithm2e}

\algref{BMTF} yields the bandlimited modified transfer function at $(\theta_n,\phi_m)$, $0 \leq n < N_\theta$, $0 \leq m < N_\phi/2$ which can be unwrapped to the remaining points by using the spherical symmetry \eqnref*{NumSphereSym}. Note that this calculation can also be performed in real-space. It is equivalent to making a Fourier interpolation of $T_k$ from $2\ell+1$ points to $N_\theta + 2\ell - 1$ points, multiplying by a bandlimited $\abs{\sin(\theta)}$, and performing a Fourier anterpolation back to $N_\theta$ points, as shown in \figref{SarvasDiag}.

\tikzset{
    >=stealth',
    textbox/.style={
           rectangle,
           draw=black,
           text centered},
    textfloat/.style={
           text centered},
}

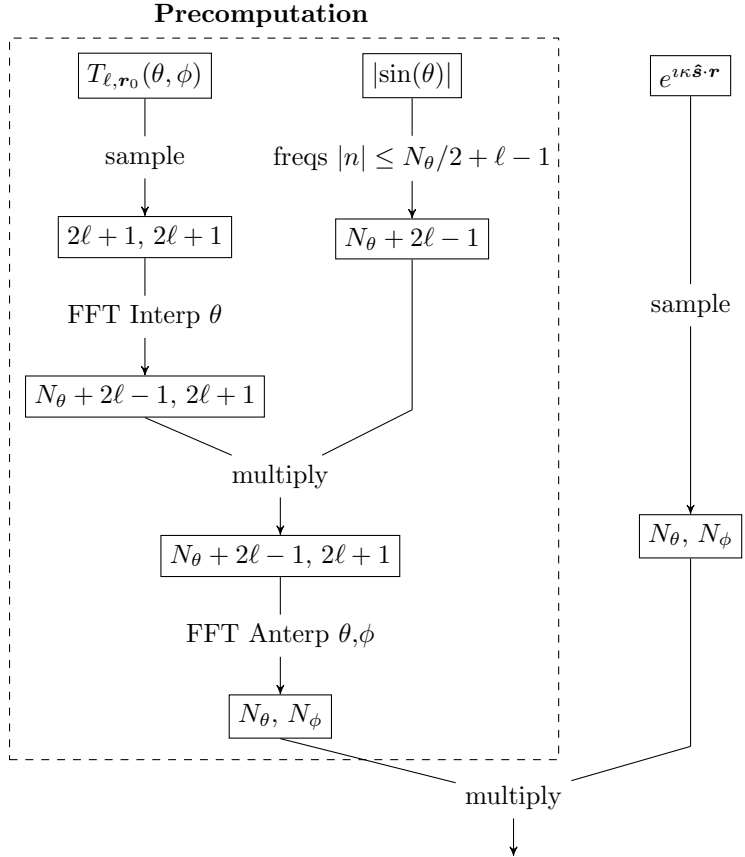
\begin{figure}[h]
\centering
\begin{tikzpicture}[node distance=1cm,auto]
\node[textbox] (T) {$T_{\ell,\v{r}_0}(\theta,\phi)$};
\node[textfloat,below=0.5cm of T] (Tsample) {sample} edge (T.south);
\node[textbox,below=0.5cm of Tsample] (ell) {$2\ell+1,\, 2\ell+1$} edge[<-] (Tsample.south);
\node[textfloat,below=0.5cm of ell] (interp) {FFT Interp $\theta$} edge (ell.south);
\node[textbox,below=0.5cm of interp] (TN) {$N_\theta + 2\ell - 1,\, 2\ell+1$} edge[<-] (interp.south);

\node[textbox,right=2cm of T] (sin) {$\abs{\sin(\theta)}$};
\node[textfloat,below=0.5cm of sin] (smooth) {freqs $\abs{n} \leq N_\theta/2 + \ell - 1$} edge (sin.south);
\node[textbox,below=0.5cm of smooth] (sinN) {$N_\theta + 2\ell - 1$} edge[<-] (smooth.south);
\node[below=2cm of sinN] (dummy) {} edge (sinN.south);

\node[textfloat,below=0.5cm of TN,xshift=1.8cm] (mult) {multiply} edge (dummy.north)
                                                                            edge (TN.south);

\node[textbox,below=0.5cm of mult] (TN2) {$N_\theta + 2\ell - 1,\, 2\ell+1$} edge[<-] (mult.south);
\node[textfloat,below=0.5cm of TN2] (anterp) {FFT Anterp $\theta$,$\phi$} edge (TN2.south);
\node[textbox,below=0.5cm of anterp] (TN3) {$N_\theta,\, N_\phi$} edge[<-] (anterp.south);

\node[textbox,right=2.5cm of sin] (eiksr) {$\eiksr$};
\node[textfloat,below=2.5cm of eiksr] (Esample) {sample} edge (eiksr.south);
\node[textbox,below=2.5cm of Esample] (EN) {$N_\theta,\, N_\phi$} edge[<-] (Esample.south);
\node[below=2.5cm of EN] (Edummy) {} edge (EN.south);

\node[textfloat,below=0.5cm of TN3,xshift=3.1cm] (mult2) {multiply} edge (TN3.south)
                                                       edge (Edummy.north);
\node[textfloat,below=0.5cm of mult2] {} edge[<-] (mult2.south);

\draw[dashed] (-1.8,0.5) rectangle (5.5,-9.1);
\node[textfloat,above=0.5cm of T,anchor=west] {\bf Precomputation};
\end{tikzpicture}
\caption{Procedure for precomputing the bandlimited modified transfer function and its application to an outgoing field. The boxed numbers give the numbers of quadrature points for $\theta$ and $\phi$ ($N_\theta$ and $N_\phi$) at each stage.}
\label{fig:SarvasDiag}
\end{figure}

Because sampling the transfer function at a single point is an $\O(\ell)$ operation and $N_\phi,N_\theta \in \O(\ell)$, the algorithm as presented is $\O(\ell^3)$. The computation of the transfer function at all sample points can be accelerated to $\O(\ell^2)$ as in~\cite{ErgulGurel_OptimalInterp} by taking advantage of its symmetry about the $\vhat{r}_0$ axis and using interpolation algorithms, but at the cost of introducing additional error.

\subsection{Choice of Quadrature}
\label{sec:DirectQuadrature}

The quadrature parameters can be constructively computed by determining the maximum error they incur. The error in computing the desired integral using the bandlimited modified transfer function with a finite uniform quadrature is
\begin{align}
\abs{\epsilon_I} &= \abs{\int_0^{2\pi} \int_0^{2\pi} E_{\v{r}}(\theta,\phi) \, \Ts_{\ell, \v{r}_0}(\theta,\phi) \d{\phi} \d{\theta} - \frac{4\pi^2}{N_\theta N_\phi} \sum_{n = 1}^{N_\theta} \sum_{m = 1}^{N_\phi} E_{\v{r}}(\theta_n,\phi_m) \, \Tsl_{\ell,\v{r}_0}(\theta_n,\phi_m)}
\label{eqn:EpsI}
\end{align}
where $\Tsl_{\ell,\v{r}_0}(\theta_n,\phi_m)$ is the bandlimited modified transfer function described in \secref{MTF}.

\subsubsection{Choosing $N_\theta$}
\label{sec:Choose_Ntheta}

The worst-case for $\epsilon_I$ in terms of $N_\theta$ occurs when $\v{r}$ and $\v{r}_0$ are aligned with the $z$-axis. This causes all spectral information to be contained in the $\theta$-direction and makes $\epsilon_I$ a function of $N_\theta$ only. Integrating $\phi$ out of \eqnref*{EpsI} yields
\begin{align*}
\abs{\epsilon^\theta_I} &= 2\pi \abs{\int_0^{2\pi} E_{\pm\abs{r}\vhat{z}}(\theta,0) \, \Ts_{\ell,\abs{\v{r}_0}\vhat{z}}(\theta,0) \d{\theta} - \frac{2\pi}{N_\theta} \sum_{n = 1}^{N_\theta} E_{\pm\abs{r}\vhat{z}}(\theta_n,0) \, \Tsl_{\ell,\abs{\v{r}_0}\vhat{z}}(\theta_n,0)}
\end{align*}
which is equivalent to the 1D case considered in \secref{FourierBackground}. Retrieving the Fourier coefficients of the plane wave in the case $\vhat{r} = \vhat{z}$ from \eqnref*{PlaneWaveSpec} and numerically computing exactly the low frequencies of the modified transfer function as described in \secref{MTF}, then \eqnref*{FourierError} leads to
\begin{align*}
\abs{\epsilon^\theta_I} = 4\pi^2 \abs{\sum_{\abs{n} \geq N_\theta/2} \FTs_n \i^n J_{\mp n}(\kappa\abs{\v{r}}) - \sum_{\abs{n} < N_\theta/2} \sum_{m \neq 0} \FTs_n \i^{n+mN_\theta} J_{\mp(n+mN_\theta)}(\kappa\abs{\v{r}})}
\end{align*}
where $\FTs_n = \F^\theta_n[\Ts_{\ell,\abs{\v{r}_0}\vhat{z}}(\theta,0)]$. Due to the very fast decay of the Bessel functions, we find it is sufficient to apply the triangle inequality and keep only the lowest order Bessel function terms:
\begin{align}
\abs{\epsilon^\theta_I} \leq 4\pi^2 \sum_{n = -\infty}^\infty \big| \FTs_n \big| \; \big| J_{M(N_\theta,n)}(\kappa \abs{\v{r}}) \big|
\label{eqn:EpsITheta}
\end{align}
where
\begin{align}
M(N_\theta,n) =
\begin{cases}
N_\theta - \abs{n} & \abs{n} < N_\theta/2\\
\abs{n} & \abs{n} \geq N_\theta/2
\end{cases}
\end{align}
\eqnref{EpsITheta} can be used to search for a value $N_\theta$ via \algref{Ntheta}.

\begin{algorithm2e}
Choose $N^{max}_\theta$ sufficiently larger than $2\ell + 1$\;
$T_n \leftarrow \Tsl_{\ell,\abs{\v{r}_0}\vhat{z}}(\frac{2\pi n}{N^{max}_\theta},0)$, $n = 0,\ldots,N^{max}_\theta-1$\;
$\tilde{T}_n \leftarrow \abs{\F^\theta_n[T]}$\;
$\tilde{E}_n \leftarrow \abs{J_n(\kappa\abs{\v{r}})}$\;
\For{$N_\theta$ from $2\ell$ to $N^{max}_\theta$ by $2$} {
	$\tilde{E}^*_n \leftarrow \tilde{E}_{M(N_\theta,n)}$\;
	\If{$\tilde{E}^* \cdot \tilde{T} < \epsilon/4\pi^2$} {return $N_\theta$}
}
\caption{Pseudocode to compute $N_\theta$ given $\ell$, $\abs{\v{r}_0}$, and $\abs{\v{r}}$.}
\label{alg:Ntheta}
\end{algorithm2e}

Since $N^{max}_\theta$ in \algref{Ntheta} is typically only a small constant larger than $2\ell+1$, the algorithm as presented is dominated by the computation of the $\O(\ell)$ modified transfer function values and requires $\O(\ell^2)$ operations. Important optimizations include using more advanced searching methods (such as bisection), applying the symmetries $\tilde{E}^*_m = \tilde{E}^*_{-m}$ and $\tilde{T}_m = \tilde{T}_{-m}$, and taking advantage of the very fast decay of $J_n$ to neglect very small terms in the dot product.

\subsubsection{Choosing $N_\phi$}
\label{sec:Choose_Nphi}

After determining an appropriate $N_\theta$, letting $N_\phi$ be a function of $\theta$ allows reducing the number of quadrature points without affecting the error. The worst-case for the integration error in terms of $N_\phi$ occurs when $\v{r}$ and $\v{r}_0$ are in the $xy$-plane. Without loss of generality, suppose $\vhat{r} = \vhat{x}$. Then, for a fixed $\theta_n$, the plane-wave \eqnref*{PlaneWaveSpec} can be expressed as
\begin{align}
E_{\v{r}}(\theta_n,\phi) = e^{\i \kappa \inner{\vhat{s}(\theta_n,\phi)}{\v{r}}} = \sum_{m=-\infty}^{\infty} i^m J_m(\kappa \abs{\v{r}} \sin(\theta_n)) \, e^{\i m \phi}
\label{eqn:PlaneWaveSpecX}
\end{align}
Since $J_m(\kappa \abs{\v{r}} \sin(\theta_n))$ is exponentially small when $m \gtrsim \kappa \abs{\v{r}} \sin(\theta_n)$, the series can be truncated at $N_\phi(\theta_n) \sim \kappa \abs{\v{r}} \sin(\theta_n)$ without incurring any appreciable error provided that the exact Fourier coefficients of $\Tsl_{\ell,\v{r}_0}(\theta,\phi)$ in $\phi$ are available. This is the case since $\Tsl_{\ell,\v{r}_0}$ is bandlimited in $\phi$. Additionally, letting $N_\phi$ be a function of $\theta$ requires a final step in the computation of the modified transfer function. \secref{MTF} computed the transfer function on a $N_\theta/2+1 \times N_\phi$ grid. With $N_{\phi} \to N_{\phi}(\theta_n)$, the data computed for each $\theta_n$ must be Fourier anterpolated from length $N_\phi$ to length $N_\phi(\theta_n)$.

Estimates of $N_\phi(\theta_n)$ can be developed by determining when $J_m(\kappa \abs{\v{r}} \sin(\theta_n))$ becomes exponentially small, as in the computation of the EBF in~\cite{Chew_Book}. However, we find that the EBF generated quadrature typically overestimates the sampling rate. To accurately compute $N_\phi(\theta_n)$ a similar procedure to that in \secref{Choose_Ntheta} is applied. After determining the appropriate $N_\theta$, the $T^{\s,L}_{\ell,\v{r}_0}$ can be computed. For a given $\theta$, the error is
\begin{align*}
\abs{\epsilon^\phi_I} &= \abs{\int_0^{2\pi} E_{\pm\abs{\v{r}}\vhat{x}}(\theta,\phi) \, \Tsl_{\ell, \pm\abs{\v{r}_0}\vhat{x}}(\theta,\phi) \d{\phi} - \frac{2\pi}{N_\phi(\theta)} \sum_{m = 1}^{N_\phi(\theta)} E_{\pm\abs{\v{r}}\vhat{x}}(\theta,\phi_m) \, \Tsll_{\ell,\pm\abs{\v{r}_0}\vhat{x}}(\theta,\phi_m)}
\end{align*}
where $\Tsll$ is the bandlimited modified transfer function with both the $\theta$-frequencies and $\phi$-frequencies truncated for the quadrature. That is, if
\begin{align*}
\Ts_{\ell,\v{r}_0}(\theta,\phi) = \frac{1}{2} T_{\ell,\v{r}_0}(\theta,\phi)\abs{\sin(\theta)} =  \sum_{n=-\infty}^{\infty} \sum_{m=-\ell}^{\ell} \FTs_{n,m} e^{n\theta + m\phi}
\end{align*}
then 
\begin{align*}
\Tsl_{\ell,\v{r}_0}(\theta,\phi) &= \sum_{n=-N_\theta/2+1}^{N_\theta/2-1} \sum_{m=-\ell}^{\ell} \FTs_{n,m} e^{n\theta + m\phi}
\intertext{and}
\Tsll_{\ell,\v{r}_0}(\theta,\phi) &= \sum_{n=-N_\theta/2+1}^{N_\theta/2-1} \sum_{m=-N_\phi(\theta_n)/2+1}^{N_\phi(\theta_n)/2-1} \FTs_{n,m} e^{n\theta + m\phi}
\end{align*}
We again apply the result of \secref{FourierBackground} by retrieving the Fourier coefficients of the plane wave in the case $\vhat{r} = \vhat{x}$ from \eqnref{PlaneWaveSpecX} and computing exactly the low frequencies of the modified transfer function. \eqnref{FourierError} leads to
\begin{align}
\abs{\epsilon^\phi_I} &\leq 2\pi \sum_{m=-\ell}^\ell \abs{\tilde{T}^{\s,L}_{m}(\theta)} \abs{J_{M(N_\phi(\theta),m)}(\kappa\abs{\v{r}}\sin(\theta))}
\label{eqn:EpsIPhi}
\end{align}

\eqnref{EpsIPhi} can then be used to search for a value of $N_\phi(\theta_n)$ via \algref{Nphi}.

\begin{algorithm2e}
Choose $N^{max}_\phi$ sufficiently larger than $2\ell + 1$\;
\For{$\theta_n$, $n = 0,\ldots,N_\theta/2$,} { 
	$T_m \leftarrow \Tsl_{\ell,\abs{\v{r}_0}\vhat{x}}(\theta_n, \frac{2\pi m}{2\ell+1})$, $m = 0,\ldots,2\ell$\;
    $\tilde{T}_m \leftarrow \abs{\F^\phi_m[T]}$\;
    $\tilde{E}_m \leftarrow \abs{J_m(\kappa\abs{\v{r}}\sin(\theta_n))}$\;
    \For{$N_\phi(\theta_n)$ from $4$ to $N^{max}_\phi$ by $4$} {
        $\tilde{E}^*_m \leftarrow \tilde{E}_{M(N_\phi(\theta_n),m)}$\;
        \If{$\tilde{E}^* \cdot \tilde{T} < \epsilon/4\pi^2$} {
            Save $N_\phi(\theta_n)$
        }
    }
}
\caption{Pseudocode to compute each $N_\phi(\theta_n)$ given $\ell$, $\abs{\v{r}_0}$, $\abs{\v{r}}$, and $N_\theta$.}
\label{alg:Nphi}
\end{algorithm2e}

Since $N^{max}_\phi$ is only a small constant larger than $2\ell+1$, the algorithm as presented is dominated by the computation of the modified transfer function and requires $\O(\ell^3)$ operations. Optimizations similar to those presented in \secref{Choose_Ntheta} can be applied. Using the EBF as an initial guess in the search for $N_\phi(\theta_n)$ further improves the searching speed. Additionally, only half of the $N_\phi(\theta_n)$'s may be computed due to symmetry about the $z = 0$ plane.

\subsubsection{Choosing $\abs{\v{r}}$ and $\abs{\v{r}_0}$}

The previous algorithms require representative values of $\abs{\v{r}}$ and $\abs{\v{r}_0}$ for each level of the tree. The worst-case transfer vectors, $\v{r}_0$, are those with smallest length. If $a_l$ is the box size at level $l$, then $\abs{\v{r}_0} = 2a_l$ is the smallest transfer vector length in the common one buffer-box case.

\begin{figure}
\centering
\includegraphics[width=2.5in]{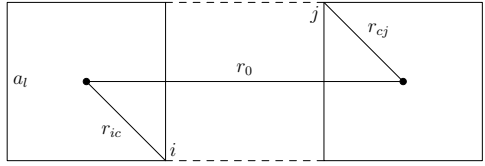}
\caption{The worst-case $\v{r}$ and $\v{r}_0$, projected from the 3D box. Here, $\abs{\v{r}_0} = 2a_l$ and $i$ and $j$ are on the opposite corners of the box so that $\abs{\v{r}} = \abs{\v{r}_{ic}} + \abs{\v{r}_{cj}} = a_l \sqrt{3}$.}
\label{fig:WorstCaseBox}
\end{figure}

The worst-case value of $\abs{\v{r}}$ is the largest. For a box of size $a_l$, $\abs{\v{r}} \leq a_l \sqrt{3}$. However, using $\abs{\v{r}} = a_l\sqrt{3}$ in the previous methods is often too conservative. This case only occurs when two points are located in the exact corners of the boxes -- a very rare case. See \figref{WorstCaseBox}. Instead, we let $\abs{\v{r}} = \alpha a_l \sqrt{3}$ for some $\alpha \in [0,1]$. A high $\alpha$ guarantees an upper bound on the error generated by the quadrature, but the points which actually generate this error become increasingly rare. A lower value of $\alpha$ will yield a smaller quadrature, but more points may fall outside the radius $\abs{\v{r}}$ for which the upper bound on the error is guaranteed.

\subsubsection{Number of Quadrature Points}
\label{sec:QuadSize}

Recall from \secref{SphereQuad} that the typical approach in the FMM is to use $N+1$ uniform points in the $\phi$ direction and $\frac{N+1}{2}$ Gauss-Legendre points in the $\theta$ direction so that all $Y^m_n$, $-n \leq m \leq n$, $0 \leq n \leq N$ are integrated exactly. In~\cite{Chew_ErrorDiagonal}, Chew et al.\ takes $\frac{N+1}{2} = \ell + 1$, which is an approximate choice based on the rapid decay of the coefficients in the spherical harmonics expansion of a plane-wave. This results in approximately
\begin{align*}
K_{sh} = 2(\ell+1)^2 \sim 2 \ell^2
\end{align*}
quadrature points.

For a given Gegenbauer series truncation $\ell$, the total number of quadrature points required in the Fourier based FMM is approximately
\begin{align*}
K_{fb} &\approx \frac{N_\theta}{2} \frac{1}{\pi} \int_{0}^{\pi}
N_\phi(\theta) \d{\theta}\\
&\approx (\ell + C_1) \frac{1}{\pi}\int_0^\pi (2\ell + C_2(\theta))\sin(\theta) \d{\theta}
\end{align*}
where $C_1,C_2 \geq 1$ are small integers dependent on $\ell$, numerically computed from the methods in \secref{Choose_Ntheta}, \secref*{Choose_Nphi}. Keeping only the leading term in $\ell$:
\begin{align*}
K_{fb} \sim \frac{4}{\pi} \,\ell^2 \approx 1.3 \, \ell^2
\end{align*}
Thus, the method presented in this paper uses approximately $0.64$ times the number of quadrature points in the standard FMM. However, it is possible that the same $N_{\phi}$ optimization can be applied to the standard FMM for the same reasons it was applied in \secref{Choose_Nphi} to reduce the standard quadrature to a comparable size.

\subsection{Interpolation and Anterpolation}
\label{sec:FFTInterp}

Most importantly, the Fourier based FMM directly uses FFTs in the interpolation and anterpolation steps. This makes the time critical upward pass and downward pass efficient and easy to implement while retaining the exactness of global methods.

\begin{figure}[htbp]
\centering
\begin{tikzpicture}[scale=0.09]
	\draw[fill=gray]  (0,21) rectangle (16,23);
	\draw[fill=gray]  (0,18) rectangle (12,26);
	\draw[fill=gray]  (0,16) rectangle (8,28);
	\draw[fill=gray]  (0,15) rectangle (4,29);
	\draw[fill=black] (0,14) rectangle (1,15);
	\draw[fill=black] (0,29) rectangle (1,30);
	\draw[step=1.0]   (0,0) grid (16,30);
	\draw (8,-3) node {m};
	\draw (-3,15) node {n};
	
	\draw[->,-stealth,thick] (18,15) -- node[anchor=south] {$\mathcal{I}^\phi$} (21,15);
	
	\draw[fill=gray]  (23+0,14) rectangle (23+16,30);
	\draw[fill=black] (23+0,14) rectangle (23+16,15);
	\draw[fill=black] (23+0,29) rectangle (23+16,30);
	\draw[step=1.0]   (23+0,0) grid (23+16,30);
	
	\draw[->,-stealth,thick] (23+18,15) -- node[anchor=south] {$\mathcal{W}$} (23+21,15);
	
	\draw[fill=gray]  (46+0,0) rectangle (46+8,30);
	\draw[fill=black] (46+0,14) rectangle (46+8,15);
	\draw[fill=black] (46+0,29) rectangle (46+8,30);
	\draw[step=1.0]   (46+0,0) grid (46+16,30);
	
	\draw[->,-stealth,thick] (46+18,15) -- node[anchor=south] {$\mathcal{A}^\theta$} (46+21,15);
	
	\draw[fill=gray]  (69+0,6) rectangle (69+8,30);
	\draw[fill=black] (69+0,17) rectangle (69+8,18);
	\draw[fill=black] (69+0,29) rectangle (69+8,30);
	\draw[step=1.0]   (69+0,0) grid (69+16,30);
	
	\draw[->,-stealth,thick] (69+18,15) -- node[anchor=south] {$\mathcal{W}$} (69+21,15);
	
	\draw[fill=gray]  (92+0,17) rectangle (92+16,30);
	\draw[fill=black] (92+0,17) rectangle (92+16,18);
	\draw[fill=black] (92+0,29) rectangle (92+16,30);
	\draw[step=1.0]   (92+0,0) grid (92+16,30);
	
	\draw[->,-stealth,thick] (92+18,15) -- node[anchor=south] {$\mathcal{A}^\phi$} (92+21,15);
	
	\draw[fill=gray]  (115+0,23) rectangle (115+12,24);
	\draw[fill=gray]  (115+0,21) rectangle (115+8,26);
	\draw[fill=gray]  (115+0,18) rectangle (115+4,29);
	\draw[fill=black] (115+0,17) rectangle (115+1,18);
	\draw[fill=black] (115+0,29) rectangle (115+1,30);
	\draw[step=1.0]   (115+0,0) grid (115+16,30);
	\draw (115+8,-3) node {m};
	\draw (115+16+3,15) node {n};
\end{tikzpicture}
\caption{The data profile at each step in an anterpolation from a large quadrature $K$ with $N_\theta = 30$ to a smaller quadrature $K'$ with $N'_\theta = 24$. The data corresponding to a pole has been darkened for clarity.}
\label{fig:InterpolateDataProfile}
\end{figure}
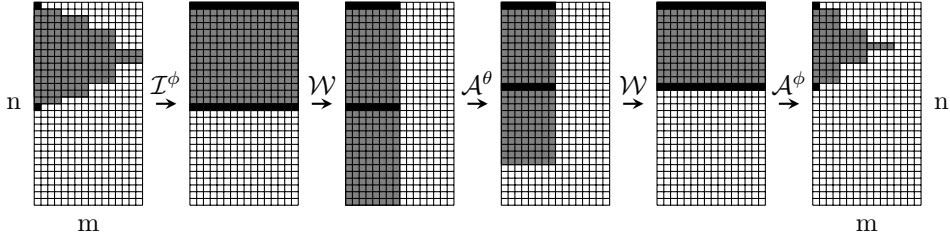

Characterize a quadrature by an array of length $N_\theta/2 + 1$,
\begin{align*}
K = [N_\phi(\theta_0), N_\phi(\theta_1), \ldots, N_\phi(\theta_{N_\theta/2-1}), N_\phi(\theta_{N_{\theta}/2})]
\end{align*}
noting that $N_\phi(\theta_n) = N_\phi(\theta_{N_\theta/2+n})$. The data $F(\theta_n,\phi_m)$ sampled on a quadrature $K$ is transformed to another quadrature $K'$ by performing a sequence of Fourier interpolations and anterpolations. Let
\begin{align*}
\mathcal{N}_\phi = \max \Big[ \; \max_{0 \leq n \leq N_\theta/2} N_\phi(\theta_n), \ \max_{0 \leq n \leq N'_\theta/2} N'_\phi(\theta_n) \; \Big]
\end{align*}
Then, the following steps, as illustrated in \figref{InterpolateDataProfile}, perform an exact interpolation/anterpolation using only FFTs.
\begin{enumerate}
\item For each $\theta_n$, $0 \leq n \leq N_\theta/2$, Fourier interpolate the data $[F(\theta_n,\phi_m)]_{m=0}^{N_\phi(\theta_n)-1}$ from length $N_\phi(\theta_n)$ to $\mathcal{N}_\phi$.

\item For each $\phi_m$, $0 \leq m < \mathcal{N}_\phi/2$, apply symmetry \eqnref*{NumSphereSym} to construct the $\theta$-periodic sequences $[F(\theta_n,\phi_m)]_{n=0}^{N_\theta-1}$.

\item For each $\phi_m$, $0 \leq m < \mathcal{N}_\phi/2$, Fourier interpolate the data $[F(\theta_n,\phi_m)]_{n=0}^{N_\theta-1}$ from length $N_\theta$ to $N'_\theta$.

\item For each $\theta_n$, $0 \leq n \leq N'_\theta/2$, apply symmetry \eqnref*{NumSphereSym} to construct the $\phi$-periodic sequences $[F(\theta_n,\phi_m)]_{m=0}^{\mathcal{N}_\phi-1}$.

\item For each $\theta_n$, $0 \leq n \leq N'_\theta/2$, Fourier anterpolate the data $[F(\theta_n,\phi_m)]_{m=0}^{\mathcal{N}_\phi-1}$ from length $\mathcal{N}_\phi$ to $N'_\phi(\theta_n)$.
\end{enumerate}
A slightly more efficient algorithm uses the symmetry \eqnref*{FourierSphereSym} rather than \eqnref*{NumSphereSym}.

\subsection{Numerical Results}
\label{sec:Results}

In the following sections, we will use the total measured error, $\epsilon_T$, defined as
\begin{align}
\epsilon_T = \frac{e^{\i\kappa\abs{\v{r}+\v{r}_0}}}{\abs{\v{r}+\v{r}_0}} - 
\sum_{n = 1}^{N_\theta} \frac{4\pi^2}{N_\theta N_\phi(\theta_n)} \sum_{m = 1}^{N_\phi(\theta_n)} E_{\v{r}}(\theta_n,\phi_m) \, \Tsl_{\ell,\v{r}_0}(\theta_n,\phi_m)
\label{eqn:eps_T}
\end{align}
The total Gegenbauer truncation error, $\epsilon_G$, is
\begin{align}
\epsilon_G = \frac{e^{\i\kappa \abs{\v{r}+\v{r}_0}}}{\abs{\v{r}+\v{r}_0}} - 
\i \kappa \sum_{n=0}^\ell (-1)^n (2n+1) h_n^{(1)}(\kappa \abs{\v{r}_0})j_n(\kappa \abs{\v{r}})P_n(\inner{\vhat{r}}{\vhat{r}_0})
\label{eqn:eps_G}
\end{align}
The total integration error $\epsilon_I$ is
\begin{align}
\epsilon_I = \epsilon_T - \epsilon_G
\label{eqn:eps_I}
\end{align}

\subsubsection{Single-Level Error}
\label{sec:SLError}

In the first test case, shown in \figref{Error_R0Z}, we examine the choice of $N_\theta$ by using the worst case $\v{r}_0 = \abs{\v{r}_0} \vhat{z}$. Here, the optimal Gegenbauer truncation $\ell$ is obtained as explained in \secref{GegenbauerTruncation} [see~\eqnref{CollinoError}]. The quadrature for computing the integral \eqnref*{FBFMM} was constructed as described in \secref{DirectQuadrature}. With box size $a = 1$, the quadrature and Gegenbauer truncation are constructed with $\abs{\v{r}} = 0.8 a \sqrt{3}$, $\abs{\v{r}_0} = 2a$, and target error $\eps = 10^{-4}$, $10^{-8}$. The plotted errors represent the maximum found over many directions $\vhat{r}$, verifying that the worst case is $\v{r} \sim \pm\vhat{z}$.

We note that the target error $\eps$ is accurately achieved, for all frequencies (except in the low-frequency breakdown because of round off errors). The actual error is within a factor of 2 or less of the target error. The increase in error for small box sizes corresponds to the low frequency breakdown when the transfer function has very large amplitude and roundoff errors become dominant. In this regime the quadrature target error bound can also be relaxed to improve efficiency --- it is inefficient to have a large quadrature that provides a small integration error when the transfer function cannot provide comparable accuracy.

We also show a comparison with a simple heuristic for choosing $N_\theta$. Based on the truncation of the Gegenbauer series, we expect $N_\theta = 2\ell + 1$ (rounded up to the nearest multiple of 2) to be a reasonable guess. The error resulting from this choice is shown by the curve labeled $\eps^{ebf}_I$. This shows that our scheme produces a more accurate estimate of this parameter and, in particular, with our scheme, we no longer see the drift in error seen with $\eps^{ebf}_I$.

\tikzset{
   eps_G/.style={blue,every mark/.append style={scale=0.7},mark=*},
   eps_I/.style={red,every mark/.append style={scale=0.7},mark=square*},
   eps_T/.style={black,every mark/.append style={fill=green!60!black},mark=x,mark size=3},
   eps_G_ebf/.style={blue!50!white,every mark/.append style={scale=0.7},mark=o},
   eps_I_ebf/.style={red!50!white,every mark/.append style={scale=0.7},mark=square},
}

\begin{figure}[ht]
\centering
\begin{tikzpicture}
    \begin{loglogaxis}[
        width = 4.25in, height = 2.5in,
        grid=both,
        legend pos = outer north east,
        xlabel=$\kappa$, xticklabel={\exptonumtickformat}, xmin=0.5, xmax=1e4,
        ytick = {1e-1, 1e-2, 1e-3, 1e-4, 1e-5, 1e-6, 1e-7, 1e-8, 1e-9, 1e-10},
        ylabel=$L_\infty$ Error, ymin=5e-10, ymax=1e-1,
    ]

	\pgfplotstableread{Data/Error_VS_Kappa_R0Z_E04.dat}\ErrorVSKappaFour

	\addplot[eps_G] table[x=kappa,y=GegenE] {\ErrorVSKappaFour};
    \addlegendentry{$\abs{\epsilon_G}$}

	\addplot[eps_I] table[x=kappa,y=IntE] {\ErrorVSKappaFour};
    \addlegendentry{$\abs{\epsilon_I}$}
    
    \addplot[eps_T] table[x=kappa,y=TotalE] {\ErrorVSKappaFour};
    \addlegendentry{$\abs{\epsilon_T}$}


    \pgfplotstableread{Data/Error_VS_Kappa_R0Z_EBF.dat}\ErrorVSKappaEBF

	\addplot[eps_I_ebf] table[x=kappa,y=IntE] {\ErrorVSKappaEBF};
    \addlegendentry{$|\epsilon_I^{ebf}|$}

    \pgfplotstableread{Data/Error_VS_Kappa_R0Z_E08.dat}\ErrorVSKappaEight

	\addplot[eps_G] table[x=kappa,y=GegenE] {\ErrorVSKappaEight};

	\addplot[eps_I] table[x=kappa,y=IntE] {\ErrorVSKappaEight};
    
    \addplot[eps_T] table[x=kappa,y=TotalE] {\ErrorVSKappaEight};
    
    
    \pgfplotstableread{Data/Error_VS_Kappa_R0Z_EBF_E08.dat}\ErrorVSKappaEBFEight

	\addplot[eps_I_ebf] table[x=kappa,y=IntE] {\ErrorVSKappaEBFEight};

    \end{loglogaxis}
\end{tikzpicture}
\caption[Fourier based Helmholtz FMM error control with $\v{r}_0 \sim \vhat{z}$]{Results with  $\v{r}_0 \sim \vhat{z}$, target accuracy $\eps = 10^{-4}$ and $\eps = 10^{-8}$, the Gegenbauer truncation chosen as in \secref{GegenbauerTruncation}, $N_\theta$ chosen as in \secref{Choose_Ntheta}, and box size $a = 1$. The $\eps^{ebf}_I$ curve corresponds to using the EBF along with $N_\theta = 2\ell + 1$; this choice leads to a drift in the error as the frequency increases. In contrast, using our error estimate and scheme to choose the parameters, the error follows closely the target error.}
\label{fig:Error_R0Z}
\end{figure}  

The second test case, shown in \figref{Error_R0X}, shows the accuracy resulting from our choice of $N_\phi$ using the worst case $\v{r}_0 = \abs{\v{r}_0}\vhat{x}$. Again, the direct computation was used to find the optimal Gegenbauer truncation $\ell$ and the quadrature was constructed following \secref{DirectQuadrature}. With box size $a = 1$, the quadrature and Gegenbauer truncation are constructed with $\abs{\v{r}} = 0.8 a \sqrt{3}$, $\abs{\v{r}_0} = 2a$, and target error $\eps = 10^{-4}$ and $10^{-8}$. The plotted errors represent the maximum found over many directions $\vhat{r}$, verifying that the worst case is $\v{r} \sim \pm\vhat{x}$. We can see that the target error is achieved even more accurately than the $\vhat{z}$ case. This is due to the number of $N_\phi$ that are chosen --- one for each $\theta_n$. The $N_\phi(\theta_n)$ that yields the most inaccurate result dominates the others, bringing the integration error very close to the target. The comparison heuristic is $\eps^{ebf}_I$ uses the constant $N_{\phi} = 2 \ell + 1$ (rounded up to the nearest multiple of 4). Again, the target error should be relaxed in the low frequency breakdown regime.

\begin{figure}[htbp]
\centering
\begin{tikzpicture}
    \begin{loglogaxis}[
        width = 4.25in, height = 2.5in,
        grid = both,
        legend pos = outer north east,
        xlabel = $\kappa$, xticklabel={\exptonumtickformat}, xmin=0.5, xmax=1e4,
        ytick = {1e-1, 1e-2, 1e-3, 1e-4, 1e-5, 1e-6, 1e-7, 1e-8, 1e-9, 1e-10},
        ylabel = $L_\infty$ Error, ymin=5e-10, ymax=1e-1,
    ]

	\pgfplotstableread{Data/Error_VS_Kappa_R0X_E04.dat}\ErrorVSKappaFour

	\addplot[eps_G] table[x=kappa,y=GegenE] {\ErrorVSKappaFour};
    \addlegendentry{$\abs{\epsilon_G}$}

	\addplot[eps_I] table[x=kappa,y=IntE] {\ErrorVSKappaFour};
    \addlegendentry{$\abs{\epsilon_I}$}
    
    \addplot[eps_T] table[x=kappa,y=TotalE] {\ErrorVSKappaFour};
    \addlegendentry{$\abs{\epsilon_T}$}
    

    \pgfplotstableread{Data/Error_VS_Kappa_R0X_EBF_E04.dat}\ErrorVSKappaEBFFour
    
    \addplot[eps_I_ebf] table[x=kappa,y=IntE] {\ErrorVSKappaEBFFour};
    \addlegendentry{$|\epsilon_I^{ebf}|$}

    \pgfplotstableread{Data/Error_VS_Kappa_R0X_E08.dat}\ErrorVSKappaEight

	\addplot[eps_G] table[x=kappa,y=GegenE] {\ErrorVSKappaEight};

	\addplot[eps_I] table[x=kappa,y=IntE] {\ErrorVSKappaEight};
    
    \addplot[eps_T] table[x=kappa,y=TotalE] {\ErrorVSKappaEight};


	\pgfplotstableread{Data/Error_VS_Kappa_R0X_EBF_E08.dat}\ErrorVSKappaEBFEight
	
	\addplot[eps_I_ebf] table[x=kappa,y=IntE] {\ErrorVSKappaEBFEight};

    \end{loglogaxis}
\end{tikzpicture}
\caption[Fourier based Helmholtz FMM error control with $\v{r}_0 \sim \vhat{x}$]{Results with $\v{r}_0 \sim \vhat{x}$, target accuracy $\eps = 10^{-4}$ and $\eps = 10^{-8}$, the Gegenbauer truncation chosen as in \secref{GegenbauerTruncation}, $N_\theta$ and $N_\phi$ chosen as in \secref{Choose_Ntheta} and \secref*{Choose_Nphi}, and box size $a = 1$. The $\eps_I^{ebf}$ curve corresponds to using the EBF along with $N_\phi(\theta_n) = 2\ell+1$; this choice greatly overestimates the size of the quadrature needed in the $\phi$-direction especially near the poles $\theta = \{0,\pi\}$. In contrast, using our error estimate and scheme to choose the parameters, the error follows closely the target error.}
\label{fig:Error_R0X}
\end{figure}

\begin{figure}[htbp]
\centering
\begin{tikzpicture}
	\begin{semilogxaxis} [
		width = 4.75in, height = 1.7in,
        grid = both, minor tick num = 1,
        xlabel = $\kappa$, xticklabel={\exptonumtickformat}, xmin=0.5, xmax=1e4,
        ylabel = $\text{QuadSize} / 2(\ell+1)^2$, ymin=0.6, ymax=1.0
	]
	
	\pgfplotstableread{Data/Error_VS_Kappa_R0X_E04.dat}\ErrorVSKappa
	
	\addplot[mark=x] table[x=kappa,y expr=\thisrow{QSize}/(2*(\thisrow{ell}+1)^2)] {\ErrorVSKappa};
	
	\end{semilogxaxis}
\end{tikzpicture}
\caption[Improved quadrature size in the Fourier based Helmholtz FMM]{The ratio of the number of quadrature points required in the Fourier based FMM and what would be used in a typical spherical harmonics based FMM for the same $\ell$. The curve asymptotes close to $2/\pi \approx 0.64$ as expected from \secref{QuadSize}}
\label{fig:QuadSize}
\end{figure}
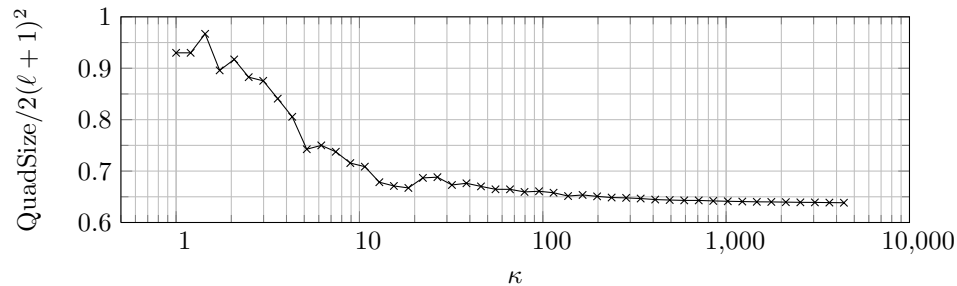

\figref{QuadSize} shows the ratio of the number of points in the quadrature used in \figref{Error_R0X} to the number of quadrature points that would be used in a typical spherical harmonics based FMM using the same Gegenbauer truncation $\ell$ chosen by the direct calculation. The procedures presented in this paper result in a quadrature which is substantially smaller than what would typically be used. Notably, the analysis in \secref{QuadSize} is supported.

Together, these results demonstrate that by choosing $\ell$ and the quadrature as presented in this paper, the error is controlled and the quadrature size is chosen nearly optimally. The accurate error bounds derived in \secref{DirectQuadrature} means that we can provide a sharp bound of the total final error of the method and optimize the running time of the method for that prescribed error. A reduction in the quadrature size improves memory usage and suggests an improved running time over similar algorithms.

\subsubsection{Multi-Level Error}
\label{sec:MLError}

The previous section verified the error bounds derived for a single level. In this section, we show that the method performs as expected in the multilevel case as well. We considered two points distributed as in \figref{WorstCaseBox}, at opposite corners of the box. The transfer pass (M2L) is always done at the highest level in the tree (level 2). Then, the number of levels is increased, thereby adding additional translation steps to the calculation. This test therefore checks that the Fourier based interpolation and anterpolation procedure (\secref{FFTInterp}) does not affect the accuracy of the calculation.

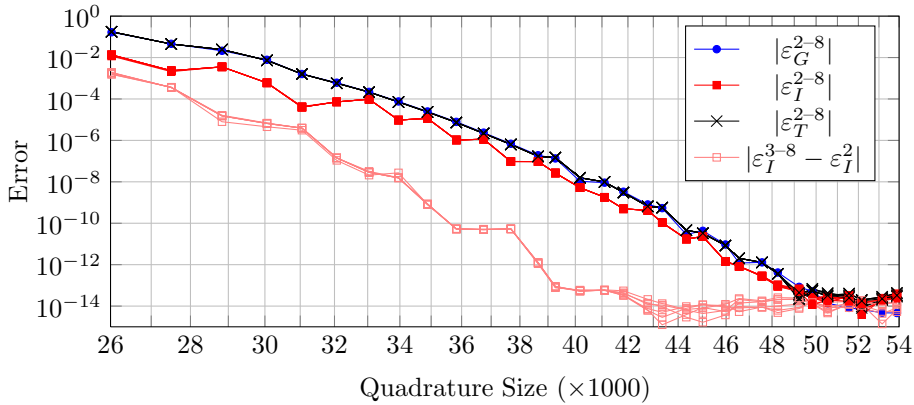
\begin{figure}[htbp]
\centering
\begin{tikzpicture}
	\begin{loglogaxis} [
		width = 4.75in, height = 2.25in,
		grid = both,
		legend pos = north east,
		xlabel = {Quadrature Size ($\times 1000$)}, 
		xtick = {26000, 27000, 28000, 29000, 30000, 31000, 32000, 33000, 34000, 35000, 36000, 37000, 38000, 39000, 40000, 41000, 42000, 43000, 44000, 45000, 46000, 47000, 48000, 49000, 50000, 51000, 52000, 53000, 54000},
		xticklabels = {26, , 28, , 30, , 32, , 34, , 36, , 38, , 40, , 42, , 44, , 46, , 48, , 50, , 52, , 54}, 
		xmin=26000, xmax=54000,
		ytick = {1e0, 1e-2, 1e-4, 1e-6, 1e-8, 1e-10, 1e-12, 1e-14},
		ylabel = {Error}, 
		ymin = 1e-15, ymax = 1e-0
	]
	
	\pgfplotstableread{Data/Error_VS_Quadrature_8.dat}\ErrorVSQuad
	
	\addplot[eps_G] table[x=QSize,y=eps_G] {\ErrorVSQuad};
	\addlegendentry{$|\epsilon_G^{2\text{--}8}|$}
	\addplot[eps_I] table[x=QSize,y=eps_I] {\ErrorVSQuad};
	\addlegendentry{$|\epsilon_I^{2\text{--}8}|$} 
	\addplot[eps_T] table[x=QSize,y=eps_T] {\ErrorVSQuad};
	\addlegendentry{$|\epsilon_T^{2\text{--}8}|$}  
	\addplot[eps_I_ebf] table[x=QSize,y=eps_A] {\ErrorVSQuad};
	\addlegendentry{$|\eps_I^{3\text{--}8} - \eps_I^2|$}  	
	
	\pgfplotstableread{Data/Error_VS_Quadrature_7.dat}\ErrorVSQuad
	\addplot[eps_I] table[x=QSize,y=eps_I] {\ErrorVSQuad};
	\addplot[eps_T] table[x=QSize,y=eps_T] {\ErrorVSQuad};
	\addplot[eps_I_ebf] table[x=QSize,y=eps_A] {\ErrorVSQuad}; 
	
	\pgfplotstableread{Data/Error_VS_Quadrature_6.dat}\ErrorVSQuad
	\addplot[eps_I] table[x=QSize,y=eps_I] {\ErrorVSQuad};
	\addplot[eps_T] table[x=QSize,y=eps_T] {\ErrorVSQuad};
	\addplot[eps_I_ebf] table[x=QSize,y=eps_A] {\ErrorVSQuad};
	
	\pgfplotstableread{Data/Error_VS_Quadrature_5.dat}\ErrorVSQuad
	\addplot[eps_I] table[x=QSize,y=eps_I] {\ErrorVSQuad};
	\addplot[eps_T] table[x=QSize,y=eps_T] {\ErrorVSQuad};
	\addplot[eps_I_ebf] table[x=QSize,y=eps_A] {\ErrorVSQuad};
	
	\pgfplotstableread{Data/Error_VS_Quadrature_4.dat}\ErrorVSQuad
	\addplot[eps_I] table[x=QSize,y=eps_I] {\ErrorVSQuad};
	\addplot[eps_T] table[x=QSize,y=eps_T] {\ErrorVSQuad};
	\addplot[eps_I_ebf] table[x=QSize,y=eps_A] {\ErrorVSQuad}; 
	
	\pgfplotstableread{Data/Error_VS_Quadrature_3.dat}\ErrorVSQuad
	\addplot[eps_I] table[x=QSize,y=eps_I] {\ErrorVSQuad};
	\addplot[eps_T] table[x=QSize,y=eps_T] {\ErrorVSQuad};
	\addplot[eps_I_ebf] table[x=QSize,y=eps_A] {\ErrorVSQuad};  
	
	\pgfplotstableread{Data/Error_VS_Quadrature_2.dat}\ErrorVSQuad
	\addplot[eps_I] table[x=QSize,y=eps_I] {\ErrorVSQuad};
	\addplot[eps_T] table[x=QSize,y=eps_T] {\ErrorVSQuad};
	\addplot[eps_I_ebf] table[x=QSize,y=eps_A] {\ErrorVSQuad};  

	\end{loglogaxis}
\end{tikzpicture}
\caption[Two-level error convergence in the Fourier based Helmholtz FMM]{Log-log plot of two-level convergence with $\kappa = 100$, $\alpha = 0.8$, and highest active level box size $a_2 = 1$. For each curve, we have 7 cases with $L=2$, \ldots, $8$. This is denoted by $2$--$8$ in the legend of the figure. The $7$ cases are nearly undistinguishable except when the error is of order $10^{-12}$ and below.}
\label{fig:Error_Quad}
\end{figure}

\figref{Error_Quad} shows the convergence of the method as the target error is adjusted, thereby changing $\ell$ and the quadrature size appropriately. We show $\eps_G^L$ (Gegenbauer truncation), $\eps_T^L$ (total error), and $\eps_I^L$ (interpolation error obtained as the difference between $\eps_G^L$ and $\eps_T^L$), for $L = 2$, \ldots, 8 ($L$ is the total number of levels). The quadrature size is the total number of quadrature points at the highest active level, with box size $a_2 = 1$. Note that for $L=2$ there is only one active level and no translation step. On the same plot, we show the discrepancy between $\eps_I^L$ ($L=3$, \ldots, 8) and $\eps_I^2$ which corresponds to the error due to the upward and downward passes. We expect this error to be much smaller than the target error, since the target error accounts for the M2L operation and therefore the quadrature is over-estimated when considering only the translation functions. This is confirmed by \figref{Error_Quad}. The curve $\eps_G^L$ is exactly the same for all $L$ while $\eps_T^L$ has small variations with $L$ due to the sampling of the translation function in the updaward and downward passes.

We note that the method converges super-exponentially, following the rate of decay of the Jacobi-Anger series \eqnref*{PlaneWaveSpec}.

\begin{figure}[htbp]
\centering
\begin{tikzpicture}
	\begin{loglogaxis} [
		width = 4.75in, height = 2.25in,
		grid = both,
		legend pos = south west,
		xlabel = {Quadrature Size ($\times 1000$)}, 
		xtick = {14000, 15000, 16000, 17000, 18000, 19000, 20000, 21000, 22000, 23000, 24000, 25000, 26000, 27000, 28000},
		xticklabels = {14, 15, 16, 17, 18, 19, 20, 21, 22, 23, 24, 25, 26, 27, 28}, 
		xmin=14000, xmax=28000, 
		ytick = {1e0, 1e-2, 1e-4, 1e-6, 1e-8, 1e-10, 1e-12, 1e-14},
		ylabel = {Error}, 
		ymin = 1e-13, ymax = 5e-2
	]
	
	\pgfplotstableread{Data/GegError4_VS_Quadrature.dat}\GegFourVSQuad
	\addplot[color=blue,mark=none,thick] table[x=QSize,y=eps_G] {\GegFourVSQuad};
	\addlegendentry{$\abs{\epsilon_G}$}	
	\addplot[eps_I] table[x=QSize,y=eps_I] {\GegFourVSQuad};
	\addlegendentry{$\abs{\epsilon_I}$}
	\addplot[eps_T] table[x=QSize,y=eps_T] {\GegFourVSQuad};
	\addlegendentry{$\abs{\epsilon_T}$}

	\pgfplotstableread{Data/GegError6_VS_Quadrature.dat}\GegSixVSQuad
	\addplot[color=blue,mark=none,thick] table[x=QSize,y=eps_G] {\GegSixVSQuad};
	\addplot[eps_I] table[x=QSize,y=eps_I] {\GegSixVSQuad};
	\addplot[eps_T] table[x=QSize,y=eps_T] {\GegSixVSQuad};
	
	\pgfplotstableread{Data/GegError8_VS_Quadrature.dat}\GegEightVSQuad
	\addplot[color=blue,mark=none,thick] table[x=QSize,y=eps_G] {\GegEightVSQuad};
	\addplot[eps_I] table[x=QSize,y=eps_I] {\GegEightVSQuad};
	\addplot[eps_T] table[x=QSize,y=eps_T] {\GegEightVSQuad};
	
	\pgfplotstableread{Data/GegError10_VS_Quadrature.dat}\GegTenVSQuad
	\addplot[color=blue,mark=none,thick] table[x=QSize,y=eps_G] {\GegTenVSQuad};
	\addplot[eps_I] table[x=QSize,y=eps_I] {\GegTenVSQuad};
	\addplot[eps_T] table[x=QSize,y=eps_T] {\GegTenVSQuad};
	
	\end{loglogaxis}
\end{tikzpicture}
\caption[Two-level error convergence to the Gegenbauer error]{Log-log plot of two-level convergence to the truncated Gegenbauer series value with $\kappa = 100$, $\alpha = 1/\sqrt{3}$, and highest active level box size $a_2 = 1$. The convergence of the numerical quadrature for $(\theta,\phi)$ (red curves) should be compared with the convergence of the toy problem in \figref{SINconv} (blue curve in \figref{SINconv}).}
\label{fig:GegError_Quad}
\end{figure}
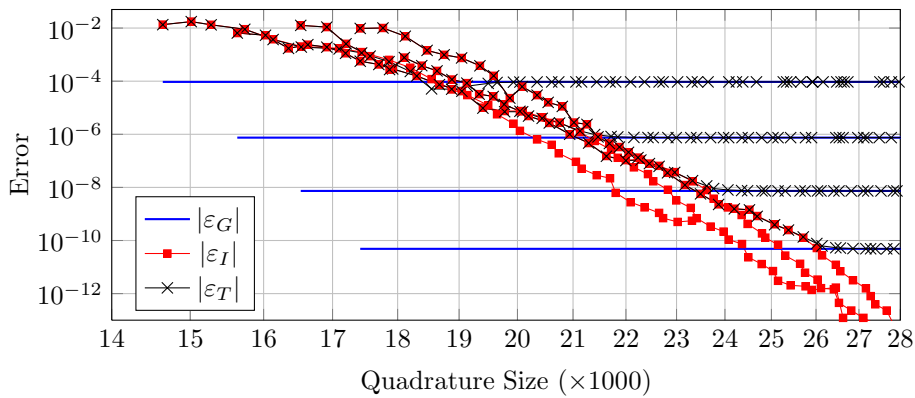

To further show that the numerical integral is converging to the Gegenbauer series value, \figref{GegError_Quad} sets the Gegenbauer error to a constant and uses the target error to increase the quadrature size only. We observe that as the quadrature size increases the total error becomes exactly equal to the Gegenbauer series truncation error. For this case, we used two active levels for the MLFMM, therefore including the FFT interpolation and anterpolation stages. This validates numerically our theoretical analysis and shows that the error caused by our anterpolation strategy and smoothing of the transfer function $\Ts_{\ell,\v{r}_0}(\theta,\phi)$ can be effectively controlled and behaves as expected.

\subsubsection{Speed}
\label{sec:Speed}

As discussed in \secref{FFTInterp}, the Fourier based FMM uses only FFTs in the upward pass and downward pass to perform the interpolations and anterpolations. FFTs make these steps easier to implement and very efficient. In addition, \secref{QuadSize} and \figref{QuadSize} suggest that the quadrature used in this paper is significantly smaller than what is typically used for a spherical harmonic basis. 

All times reported are from repeated runs on a 2.93GHz Intel Core i7 Quad Core Processor 870 with 8MB Cache with 8GB of 1333MHz DDR3 RAM.

First, we illustrate the efficient use of FFTs in the interpolation and anterpolation stages. Here, we compare the FFT interpolation method of \secref{FFTInterp} against the semi-naive exact spherical harmonic interpolation method which consists of a forward FFT in $\phi$, a dense matrix-matrix product on the $\theta$ angles, and a backward FFT in $\phi$. This spherical harmonics method is analyzed and accelerated in \cite{AlpertJokobChien_SphereFilter,Darve_Numerical}. \figref{InterpTime} confirms the expected asymptotic running times of each method.

\begin{figure}[ht]
\centering
\begin{tikzpicture}	
	\begin{loglogaxis} [
		width = 4.75in, height = 2.5in,
        grid = both,
        /pgf/number format/fixed,
        xlabel = {Gegenbauer Truncation, $\ell$},
        xtick={50,100,200,500,1000,2000,5000},
        xticklabel={\exptonumtickformat}, xticklabel style={/pgf/number format/.cd, fixed, fixed zerofill, precision=0}, 
        xmin=50, xmax=6000,
        ylabel = {CPU Time (sec)}, yticklabel={\exptonumtickformat}, yticklabel style={/pgf/number format/.cd, fixed, precision=5}, 
        ymin = 2e-4, ymax = 5e3,
        ytick={1e-5,1e-4,1e-3,1e-2,1e-1,1e0,1e1,1e2,1e3},
		legend pos=north west
	]
	
	
	

	\pgfplotstableread{Data/InterpTimes3.dat}\Interp
	
	\addplot[mark=square*,thick,blue] table[x=ell,y=PolyInterp] {\Interp};
	\addlegendentry{Semi-naive}
	
	\addplot[mark=*,thick,red] table[x=ell,y=FFTInterp3] {\Interp};
	\addlegendentry{FFT}

	\addplot[mark=none,dashed,very thick] coordinates { (1000,1) (1e4,1000) };
	\addlegendentry{\scriptsize $\O(\ell^3)$}
	\addplot[mark=none,dotted,very thick] coordinates { (1000,1) (1e4,230) };
	\addlegendentry{\scriptsize $\O(\ell^2 \log(\ell))$}

	\end{loglogaxis}
	
\end{tikzpicture}
\caption{Comparison of the FFT interpolation scheme of \secref{FFTInterp} [complexity $\O(\ell^2 \log \ell)$] with the semi-naive $\O(\ell^3)$ spherical harmonic transform described in~\cite{AlpertJokobChien_SphereFilter}.
}
\label{fig:InterpTime}
\end{figure}

\begin{figure}[ht]
\centering
\begin{tikzpicture}	
	\begin{loglogaxis} [
		width = 4.75in, height = 3in,
        grid = both,
        /pgf/number format/fixed,
        xlabel = {(N,MHz)}, 
        xmin=1e3, xmax=1e7,
        xtick={1e3,1e4,1e5,1e6,1e7,1e8},
        xticklabels={(1K,23.9),(10K,51.4),(100K,111),(1M,239),(10M,514)},
        ylabel = {CPU Time (sec)}, yticklabel={\exptonumtickformat}, ymin = 0.03, ymax = 1e5,
		legend pos=north west
	]
	
	\pgfplotstableread{Data/NLogN.dat}\NLogN
	
	\addplot[mark=none,very thick] table[x=N,y=DIRtime] {\NLogN};
	\addlegendentry{\small $\O(N^2)$ Direct}
	\addplot[mark=none,dashed,very thick] coordinates { (1e3,0.05) (1e8,5000) };
	\addlegendentry{\small $\O(N)$}
	
	\addplot[mark=x,blue] table[x=N,y=FMMtime2] {\NLogN};
	\addlegendentry{\small 1 Level}
	\addplot[mark=triangle*,red] table[x=N,y=FMMtime3] {\NLogN};
	\addlegendentry{\small 2 Levels}
	\addplot[mark=square*,green!60!black] table[x=N,y=FMMtime4] {\NLogN};
	\addlegendentry{\small 3 Levels}
	\addplot[mark=pentagon*,purple] table[x=N,y=FMMtime5] {\NLogN};
	\addlegendentry{\small 4 Levels}
	\addplot[mark=*,orange] table[x=N,y=FMMtime6] {\NLogN};
	\addlegendentry{\small 5 Levels}

	\end{loglogaxis}
	
\end{tikzpicture}
\caption[$\O(N)$ running time of the Fourier based Helmholtz FMM]{Average running times of the Fourier based FMM for constant number of volumetric points per wavelength. By choosing the optimal number of levels, we achieve an $\O(N)$ complexity.}
\label{fig:RunningTime}
\end{figure}
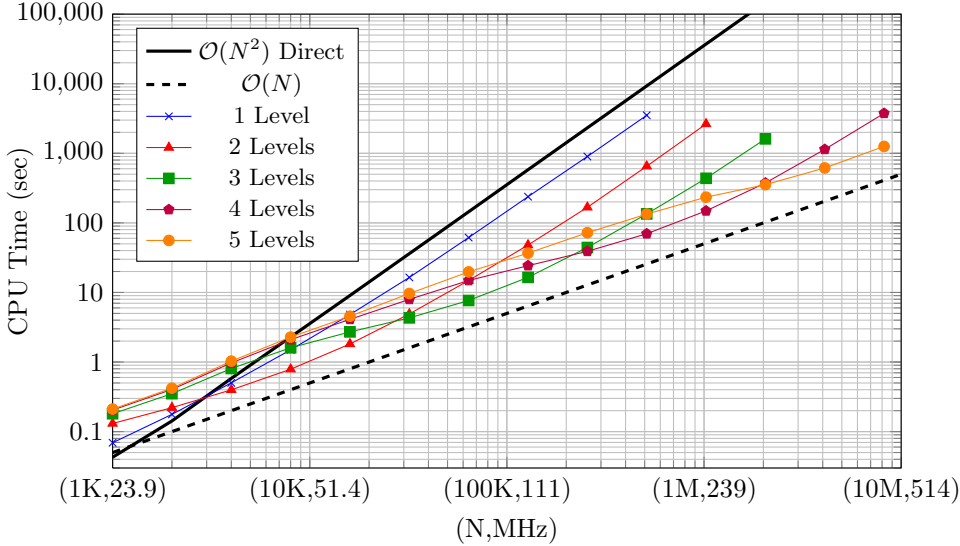

To show that the optimal asymptotic running time is achieved, \figref{RunningTime} shows the recorded running times of the Fourier based FMM and the direct matrix-vector product. The target error is set to $10^{-4}$ with $\alpha = 1$ and is achieved in every case. For $N=8.2 \cdot 10^6$ the points are uniformly distributed in a cube with side length $80\lambda$. The wave number $\kappa$ is scaled with $N^{1/3}$ to provides a nearly constant density of points per wavelength as $N$ is varied. As expected, by choosing the correct number of levels the running time is asymptotically $\O(N)$ as $N$ is increased with a constant number of points per wavelength. Note that the cross-over point is $N \approx 3,000$. 

\section{Conclusion}

We have proposed using the Fourier basis $e^{\i p \phi}e^{\i q \theta}$ in the spherical variables $\phi$ and $\theta$ to represent the far-field approximation in the FMM. By approximating the Helmholtz kernel with \eqnref{modifiedIntegrand}
and using a uniform quadrature we can take advantage of very fast, exact, and well-known FFT interpolation/anterpolation methods. By exploiting symmetries and a scheme to reduce the number of points in the $\phi$-direction, the total number of uniform quadrature points required is smaller than the number of Gauss-Legendre quadrature points typically used with spherical harmonics. This is realized by removing the high frequency components of the modified transfer function, $\Ts_{\ell,\v{r}_0}(\theta,\phi)$, during the precomputation phase which do not significantly contribute to the final integration.

The Fourier based FMM approach has a number of advantages. Since the interpolation and anterpolation algorithms are exact, the error analysis is simplified; we establish a sharp upper bound for the error. The key parameters are the Gegenbauer truncation parameter $\ell$ and the quadrature size, in particular the sampling rate in the $\theta$-direction. The truncation error $\epsilon_G$ has been extensively studied by other authors and is well understood. The integration error $\epsilon_I$ accounts for the bandlimited approximation of the modified transfer function and the finite sampling of the plane-waves. This error can be accounted for a priori during the precomputation stage. Numerical tests have confirmed that this error analysis is quite sharp. Constructive algorithms to find nearly optimal parameters were proposed. Since highly efficient FFT algorithms are available in virtually every computing environment, the time-critical interpolation stages of the algorithm are much easier to implement efficiently.

\section*{Acknowledgments}
\label{sec:acknowledgements}

This research was supported by the U.S.\ Army Research Laboratory, through the Army High Performance Computing Research Center, Cooperative Agreement W911NF-07-0027, the Stanford School of Engineering, and the King Abdullah University of Science and Technology.

\bibliographystyle{siam}
\bibliography{bibfile}

\begin{thebibliography}{10}

\bibitem{Collino_JA}
{\sc Q.~Carayol and F.~Collino}, {\em Error estimates in the fast multipole
  method for scattering problems. {Part} 1: Truncation of the {Jacobi-Anger}
  series}, ESAIM: M2NA, 38 (2004), pp.~371--394.

\bibitem{Collino_GS}
\leavevmode\vrule height 2pt depth -1.6pt width 23pt, {\em Error estimates in
  the fast multipole method for scattering problems. {Part} 2: Truncation of
  the {Gegenbauer} series}, ESAIM: M2NA, 39 (2004), pp.~183--221.

\bibitem{Chew_Book}
{\sc W.C. Chew, E.~Michielssen, J.~M. Song, and J.~M. Jin}, eds., {\em Fast and
  Efficient Algorithms in Computational Electromagnetics}, Artech House, Inc.,
  Norwood, MA, USA, 2001.

\bibitem{Chowdhury_FST}
{\sc Indranil Chowdhury and Vikram Jandhyala}, {\em Integration and
  interpolation based on fast spherical transforms for the multilevel fast
  multipole method}, Microwave and Optical Technology Letters, 48 (2006),
  pp.~1961--1964.

\bibitem{Rokhlin_FMMPedestrian}
{\sc R.~Coifman, V.~Rokhlin, and S.~Wandzura}, {\em The fast multipole method
  for the wave equation: A pedestrian prescription}, Antennas and Propag.
  Magazine, IEEE, 35 (1993), pp.~7--12.

\bibitem{Darve_ErrorAsymptotic}
{\sc Eric Darve}, {\em The fast multipole method {I}: Error analysis and
  asymptotic complexity}, SIAM J. Numer. Anal., 38 (2000), pp.~98--128.

\bibitem{Darve_Numerical}
\leavevmode\vrule height 2pt depth -1.6pt width 23pt, {\em The fast multipole
  method: Numerical implementation}, J. Comput. Phys., 160 (2000),
  pp.~195--240.

\bibitem{DriscollHealy_FST}
{\sc James~R. Driscoll and D.M. Healy}, {\em Computing {Fourier} transforms and
  convolutions on the 2-sphere}, Adv. Appl. Math., 15 (1994), pp.~202--250.

\bibitem{Rokhlin_FastPolynomialInterp}
{\sc A.~Dutt, M.~Gu, and V.~Rokhlin}, {\em Fast algorithms for polynomial
  interpolation, integration, and differentiation}, SIAM J. Numer. Anal., 33
  (1996), pp.~1689--1711.

\bibitem{Rokhlin_FMM_EM}
{\sc N.~Engheta, W.~D. Murphy, V.~Rokhlin, and M.~S. Vassiliou}, {\em The fast
  multipole method {(FMM)} for electromagnetic scattering problems}, IEEE
  Trans. Antennas Propag., 40 (1992), p.~634.

\bibitem{ErgulGurel_OptimalInterp}
{\sc O.~Ergul and L.~Gurel}, {\em Optimal interpolation of translation operator
  in multilevel fast multipole algorithm}, IEEE Trans. Antennas Propag., 54
  (2006), pp.~3822--3826.

\bibitem{Chew_ErrorControl}
{\sc M.~L. Hastriter, S.~Ohnuki, and W.~C. Chew}, {\em Error control of the
  translation operator in {3D MLFMA}}, Microwave and Optical Technology
  Letters, 37 (2003), pp.~184--188.

\bibitem{HealyMoore_FST}
{\sc D.M. Healy, D.~Rockmore, P.J. Kostelec, and S.~Moore}, {\em {FFTs for the
  2-sphere} - improvements and variations}, J. Fourier Analysis and
  Applications, 9 (2003), pp.~341--385.

\bibitem{AlpertJokobChien_SphereFilter}
{\sc R\"{u}diger Jakob-Chien and Bradley~K. Alpert}, {\em A fast spherical
  filter with uniform resolution}, J. Comput. Phys., 136 (1997), pp.~580--584.

\bibitem{Knab_InterpProlateSpheroid}
{\sc J.~Knab}, {\em Interpolation of bandlimited functions using the
  approximate prolate series (corresp.)}, Information Theory, IEEE, 25 (1979),
  pp.~717--720.

\bibitem{Chew_ErrorDiagonal}
{\sc S.~Koc, J.~M. Song, and W.~C. Chew}, {\em Error analysis for the numerical
  evaluation of the diagonal forms of the scalar spherical addition theorem},
  SIAM J. Numer. Anal., 36 (1999), pp.~906--921.

\bibitem{mckay2005fast}
{\sc E.~McKay~Hyde and OP~Bruno}, {\em {A fast higher-order solver for
  scattering by penetrable bodies in three dimension}}, J. Comp. Phys., 202
  (2005), pp.~236--261.

\bibitem{Rahola_DiagonalForms}
{\sc J.~Rahola}, {\em Diagonal forms of the translation operators in the fast
  multipole algorithm for scattering problems}, BIT, 36 (1996), pp.~333--358.

\bibitem{Rokhlin_FMMScattering}
{\sc V.~Rokhlin}, {\em Rapid solution of integral equations of scattering
  theory in two dimensions}, J. Comput. Phys., 86 (1990), pp.~414--439.

\bibitem{Rokhlin_DiagonalTranslation}
\leavevmode\vrule height 2pt depth -1.6pt width 23pt, {\em Diagonal forms of
  translation operators for the {Helmholtz} equation in three dimensions},
  Applied and Computational Harmonic Analysis, 1 (1993), pp.~82--93.

\bibitem{RokhlinTygert_FST}
{\sc V.~Rokhlin and M.~Tygert}, {\em Fast algorithms for spherical harmonic
  expansions}, SIAM J. Scientific Computing, 27 (2006), pp.~1903--1928.

\bibitem{Sarvas_FFT}
{\sc Jukka Sarvas}, {\em Performing interpolation and anterpolation entirely by
  fast {Fourier} transform in the {3-D} multilevel fast multipole algorithm},
  SIAM J. Numer. Anal., 41 (2003), pp.~2180--2196.

\bibitem{Sneeuw_2DFourierSphere}
{\sc Nico Sneeuw and Richard Bun}, {\em Global spherical harmonic computation
  by two-dimensional fourier methods}, Journal of Geodesy, 70 (1996),
  pp.~224--232.
\newblock 10.1007/BF00873703.

\bibitem{Suda_FST}
{\sc Reiji Suda and Masayasu Takami}, {\em A fast spherical harmonics transform
  algorithm}, Math. of Comp., 71 (2001), pp.~703--715.

\end{thebibliography}

\end{document}